\documentclass[12pt]{article}
\usepackage[utf8x]{inputenc}
\usepackage{amssymb, amsmath, amsthm, amscd}
\usepackage[dvips]{graphics}
\usepackage[all]{xy}
\usepackage{bbm}
\usepackage{enumitem}
\usepackage{setspace}
\usepackage[colorlinks=true]{hyperref}
\hypersetup{colorlinks   = true}
\hypersetup{linkcolor=blue}
\begin{document}
\newtheorem{The}{Theorem}[section]
\newtheorem{Lem}[The]{Lemma}
\newtheorem{Prop}[The]{Proposition}
\newtheorem{Cor}[The]{Corollary}
\newtheorem{Rem}[The]{Remark}
\newtheorem{Obs}[The]{Observation}
\newtheorem{SConj}[The]{Standard Conjecture}
\newtheorem{Titre}[The]{\!\!\!\! }
\newtheorem{Conj}[The]{Conjecture}
\newtheorem{Question}[The]{Question}
\newtheorem{Prob}[The]{Problem}
\newtheorem{Def}[The]{Definition}
\newtheorem{Not}[The]{Notation}
\newtheorem{Claim}[The]{Claim}
\newtheorem{Conc}[The]{Conclusion}
\newtheorem{Ex}[The]{Example}
\newtheorem{Fact}[The]{Fact}
\newtheorem{Formula}[The]{Formula}
\newtheorem{Formulae}[The]{Formulae}
\newcommand{\C}{\mathbb{C}}
\newcommand{\R}{\mathbb{R}}
\newcommand{\N}{\mathbb{N}}
\newcommand{\Z}{\mathbb{Z}}
\newcommand{\Q}{\mathbb{Q}}
\newcommand{\Proj}{\mathbb{P}}
\newcommand{\Rc}{\mathcal{R}}
\newcommand{\Oc}{\mathcal{O}}

\begin{center}

{\Large\bf Adiabatic Limit and the Fr\"olicher Spectral Sequence}

\end{center}

\begin{center}

{\large Dan Popovici}

\end{center}

\vspace{1ex}

\noindent{\small{\bf Abstract.} Motivated by our conjecture of an earlier work predicting the degeneration at the second page of the Fr\"olicher spectral sequence of any compact complex manifold supporting an SKT metric $\omega$ (i.e. such that $\partial\bar\partial\omega=0$), we prove degeneration at $E_2$ whenever the manifold admits a Hermitian metric whose torsion operator $\tau$ and its adjoint vanish on $\Delta''$-harmonic forms of positive degrees up to $\mbox{dim}_\C X$. Besides the pseudo-differential Laplacian inducing a Hodge theory for $E_2$ that we constructed in earlier work and Demailly's Bochner-Kodaira-Nakano formula for Hermitian metrics, a key ingredient is a general formula for the dimensions of the vector spaces featuring in the Fr\"olicher spectral sequence in terms of the asymptotics, as a positive constant $h$ decreases to zero, of the small eigenvalues of a rescaled Laplacian $\Delta_h$, introduced here in the present form, that we adapt to the context of a complex structure from the well-known construction of the adiabatic limit and from the analogous result for Riemannian foliations of \'Alvarez L\'opez and Kordyukov.}

\vspace{1ex}

\section{Introduction}\label{section:introd}

Let $X$ be a compact complex manifold of dimension $n$. It is well known that the existence of a K\"ahler metric $\omega$ on $X$ implies the degeneration at $E_1$ of the Fr\"olicher spectral sequence that relates the complex structure of $X$ (encapsulated in the Dolbeault, i.e. the $\bar\partial$-, cohomology $H^{p,,q}(X,\,\C)$, the start page of this spectral sequence) to the differential structure of $X$ (encapsulated in the De Rham, i.e. the $d$-, cohomology $H^k(X,\,\C)$, the limiting page of this spectral sequence). However, since K\"ahler metrics exist only rarely when $n\geq 3$, it is natural to search for weaker metric conditions on $X$ that ensure a (possibly weaker) degeneration property of the algebro-geometric object that is the Fr\"olicher spectral sequence of $X$. The best we can hope for in the non-K\"ahler context is the degeneration at the second page. To this end, we proposed the following conjecture in [Pop16]:

\begin{Conj}\label{Conj:SKT_E2} If a compact complex manifold $X$ admits an {\bf SKT metric} $\omega$ (i.e. a Hermitian metric $\omega$ such that $\partial\bar\partial\omega=0$), the Fr\"olicher spectral sequence of $X$ degenerates at $E_2$.

\end{Conj}

There is evidence that this ought to be true. The statement holds true on all the examples of compact complex manifolds that we are aware of, namely all the $3$-dimensional nilmanifolds, the $3$-dimensional solvmanifolds that are currently classified, the Calabi-Eckmann manifold $S^3\times S^3$, etc. In [Pop16], we proved this statement under the extra assumption that the SKT metric $\omega$ which is supposed to exist has a {\it small torsion} in the sense that the upper bound of its torsion operator of type $(0,\,0)$ (defined in a precise way) does not exceed a third of the spectral gap of the elliptic, self-adjoint and non-negative, differential operator $\Delta' + \Delta''$ in every bidegree $(p,\,q)$. As usual, $\Delta'=\Delta'_\omega=\partial\partial^\star_\omega + \partial^\star_\omega\partial$ and $\Delta''=\Delta''_\omega=\bar\partial\bar\partial^\star_\omega + \bar\partial^\star_\omega\bar\partial$ are the $\partial$-, resp. $\bar\partial$-Laplacians on smooth differential forms on $X$.    

While Conjecture \ref{Conj:SKT_E2} remains elusive at the moment, we give in this paper a different sufficient metric condition for degeneration at $E_2$ that does not assume the fixed Hermitian metric $\omega$ to be SKT. As usual (see e.g. [Dem84] or [Dem97, VII, $\S.1$]), we consider the torsion operator $\tau=\tau_\omega:=[\Lambda_\omega,\,\partial\omega\wedge\cdot]$ of type $(1,\,0)$ defined on smooth differential forms on $X$, where $\Lambda_\omega$ is the adjoint of the multiplication by $\omega$ w.r.t. the inner product defined by $\omega$, while $[A,\,B]=AB-(-1)^{ab}\,BA$ is the graded commutator of any two endomorphisms $A,B$ of respective degrees $a, b$ of the bi-graded algebra $C^\infty_{\bullet,\,\bullet}$ of smooth differential forms on $X$. Specifically, we prove

\begin{The}\label{The:main2} Let $(X,\,\omega)$ be a compact Hermitian manifold with $\mbox{dim}_\C X=n$ such that the inclusion of kernels \begin{equation}\label{eqn:torsion_hyp}\ker\Delta''\subset\ker\,[\tau,\,\tau^\star]\end{equation}

\noindent holds for the operators $\Delta'', [\tau,\,\tau^\star]:C^\infty_k(X,\,\C)\longrightarrow C^\infty_k(X,\,\C)$ in every degree $k\in\{1,\dots , n\}$.

 Then, the Fr\"olicher spectral sequence of $X$ degenerates at the second page $E_2$.

\end{The}

Hypothesis (\ref{eqn:torsion_hyp}) is of a qualitative nature and it is comparatively easy to check on concrete examples of compact Hermitian manifolds $(X,\,\omega)$ whether it holds or not. For example, $S^3\times S^3$ equipped with the Calabi-Eckmann complex structure and the Iwasawa manifold do not satisfy it when they are given the natural non-K\"ahler metrics (easy verifications that are left to the reader). Intuitively, (\ref{eqn:torsion_hyp}) requires the torsion of $\omega$ to be ``small'' since, for non-negative operators, the smaller one has a larger kernel. (We will use throughout the paper the usual order relation for linear operators $A,B$: $A\geq B$ will mean that $\langle\langle Au,\,u\rangle\rangle\geq\langle\langle Bu,\,u\rangle\rangle$ for all forms $u$, where $\langle\langle\,\,,\,\,\rangle\rangle$ stands for the $L^2$ inner product induced by the fixed Hermitian metric $\omega$ on $X$.) Hypothesis (\ref{eqn:torsion_hyp}) is obviously satisfied if $\omega$ is K\"ahler since $\tau=0$ in that case. We do not know whether there exist compact complex non-K\"ahler manifolds that satisfy hypothesis (\ref{eqn:torsion_hyp}). 

\vspace{3ex}

Inspired by the extensive literature on the adiabatic limit associated with a Riemannian foliation (see e.g. [Wi85], [MM90], [For95], [ALK00] and the references therein), we adapt that construction to the case of the splitting $d=\partial+\bar\partial$ defining the complex structure of $X$. Thus, for every constant $h>0$ that is eventually let to converge to $0$, we define in section $\S.$\ref{section:rescaled-laplacians} two {\bf rescalings} of the usual $d$-Laplacian $\Delta=dd^\star + d^\star d$ acting on the smooth differential forms on an arbitrary compact Hermitian manifold $(X,\,\omega)$: 

$$\Delta_h:=d_hd_h^\star + d_h^\star d_h,$$

\noindent where $d_h:=h\partial + \bar\partial$ modifies $d$ by {\bf rescaling $\partial$} while keeping $\bar\partial$ fixed, but its formal adjoint $d_h^\star$ is computed w.r.t. the given Hermitian metric $\omega$, and 

$$\Delta_{\omega_h}:=dd_{\omega_h}^\star + d_{\omega_h}^\star d,$$

\noindent where $d=\partial + \bar\partial$ is kept unchanged, but its formal adjoint $d_{\omega_h}^\star$ is computed w.r.t. a {\bf rescaled metric $\omega_h$} that modifies the original $\omega$ by multiplying the pointwise inner product of $(p,\,q)$-forms by $h^{2p}$. So, the anti-holomorphic degree $q$ of $(p,\,q)$-forms does not contribute to the definition of $\omega_h$. Although strongly inspired by the adiabatic limit construction in the presence of a Riemannian foliation, this partial rescaling of a Hermitian metric seems to be new and to hold further promise for the future.  

 In section $\S.$\ref{section:rescaled-laplacians}, we study these two {\bf rescaled Laplacians} and the relationships between them. As in the foliated case of [ALK00], $\Delta_h$ and $\Delta_{\omega_h}$ are seen to have the same spectrum and to have eigenspaces that are obtained from each other via a {\it rescaling isometry}.

A key ingredient in the proof of Theorem \ref{The:main2} is the following formula for the dimensions of the vector spaces featuring on each page of the Fr\"olicher spectral sequence of $X$ in terms of the number of small eigenvalues of the rescaled Laplacian $\Delta_h$ (or, equivalently, $\Delta_{\omega_h}$). ``Small'' refers to the eigenvalues' decay rate to zero as $h\downarrow 0$. This result and its proof are strongly inspired by the analogous result for foliations proved by \'Alvarez L\'opez and Kordyukov in [ALK00]. However, to our knowledge, this particular form of the result in the context of the Fr\"olicher spectral sequence seems new and is of independent interest.

\begin{The}\label{The:main1} Let $(X,\,\omega)$ be a compact Hermitian manifold with $\mbox{dim}_\C X=n$. For every $r\in\N^\star$ and every $k=0,\dots , 2n$, the following identity holds:

\begin{equation}\label{eqn:main_identity}\mbox{dim}_\C E_r^k = \sharp\{i\,\mid\,\lambda_i^k(h)\in O(h^{2r}) \hspace{2ex} \mbox{as} \hspace{1ex} h\downarrow 0\},\end{equation}

\noindent where $E_r^k:=\oplus_{p+q=k}E_r^{p,\,q}$ is the direct sum of the spaces of total degree $k$ on the $r^{th}$ page of the Fr\"olicher spectral sequence of $X$, while $0\leq \lambda_1^k(h)\leq \lambda_2^k(h)\leq\dots\leq\lambda_i^k(h)\leq\dots$ are the eigenvalues, counted with multiplicities, of the rescaled Laplacian $\Delta_h : C^\infty_k(X,\,\C)\longrightarrow C^\infty_k(X,\,\C)$ ($=$ those of $\Delta_{\omega_h}: C^\infty_k(X,\,\C)\longrightarrow C^\infty_k(X,\,\C)$) acting on $k$-forms. As usual, $\sharp$ stands for the cardinal of a set.

\end{The}

The proof of this statement proceeds along the lines of the one given in [ALK00] for the analogous statement in the foliated case with some simplifications, adjustments and inevitable differences in detail. We spell it out in section $\S.$\ref{section:Laplacians_Froelicher}. In the proof of Theorem \ref{The:main1}, we also use our pseudo-differential Laplacian $\widetilde\Delta=\partial p''\partial^\star + \partial^\star p''\partial + \Delta'' :C^\infty_{p,\,q}(X,\,\C)\longrightarrow C^\infty_{p,\,q}(X,\,\C)$ (where $p''$ is the orthogonal projection onto $\ker\Delta''$) constructed in every bidegree $(p,\,q)$ in [Pop16] and shown there to induce a Hodge isomorphism between its kernel and the space $E_2^{p,\,q}$ of bidegree $(p,\,q)$ featuring on the second page of the Fr\"olicher spectral sequence.

Along with Theorem \ref{The:main1} and the pseudo-differential Laplacian $\widetilde\Delta$, the third main ingredient in the proof of Theorem \ref{The:main2} is the following formula of the Bochner-Kodaira-Nakano type for Hermitian (not necessarily K\"ahler) metrics $\omega$ established by Demailly in [Dem84] (see also [Dem97, VII, $\S.1$):

\begin{equation}\label{BKN_introd}\Delta'' = \Delta'_\tau + [\Lambda,\,[\Lambda,\,\frac{i}{2}\partial\bar\partial\omega]] - [\partial\omega\wedge\cdot,\,(\partial\omega\wedge\cdot)^\star],\end{equation}

\noindent where $[\bullet,\bullet]$ is the usual graded commutator (see e.g. Notation \ref{Not:introd} below), $\Lambda=\Lambda_\omega$ is the adjoint of the multiplication operator $\omega\wedge\cdot$, $\tau=\tau_\omega:=[\Lambda,\,\partial\omega\wedge\cdot]$ is the torsion operator of $\omega$ and $\Delta'_\tau:=[\partial+\tau,\,(\partial+\tau)^\star]$. This formula enables us to compare various Laplacians and finish the proof of Theorem \ref{The:main2} in section $\S.$\ref{section:degeneration_E2}.

This paper owes much to the ideas and techniques in our main source of inspiration [ALK00] and to the treatment given to the Leray spectral sequence in [MM90] and [For95], although the setting and the objectives are different.  

 In the Appendix, we give an estimate of the discrepancy between the Laplacians $\Delta'$ and $\Delta''$ under the SKT assumption on the metric $\omega$ (cf. Lemma \ref{Lem:Delta'-''_comparison_SKT}). This is of independent interest and leads to the lower bound $-Ch^2$ for the operator $\Delta_h-h^2\Delta$ for all $0<h<1$ when $\omega$ is SKT, where $C\geq 0$ is a constant independent of $h$ that can be chosen to be any upper bound of the non-negative {\bf bounded} torsion operator $[\bar\tau,\,\bar\tau^\star]$ (cf. Lemma \ref{Lem:Delta'-''_comparison_SKT_2}). In view of Theorem \ref{The:main1} and some minor extra arguments, if the lower bound  $-Ch^2$ could be improved to $0$, Conjecture \ref{Conj:SKT_E2} would be solved, but at the moment we are unfortunately short of arguments to perform this improvement.  

\begin{Not}\label{Not:introd}{\rm For a given Hermitian metric $\omega$ on a given compact complex manifold $X$, $\langle\langle\,\,,\,\,\rangle\rangle = \langle\langle\,\,,\,\,\rangle\rangle_\omega$ will stand for the $L^2$ inner product defined by $\omega$ on the spaces $C^\infty_{p,\,q}(X,\,\C)$ (resp. $C^\infty_k(X,\,\C)$) of smooth differential $(p,\,q)$-forms (resp. $k$-forms) on $X$, while $||\,\,||=||\,\,||_\omega$ will denote the corresponding $L^2$-norm. For self-adjoint linear operators $A,B$ on the bi-graded algebra $\oplus_{p,\,q}C^\infty_{p,\,q}(X,\,\C)$, by $A\geq B$ we shall mean (as is the standard convention) that $\langle\langle Au,\,u\rangle\rangle\geq\langle\langle Bu,\,u\rangle\rangle$ for every form $u$ lying in the space on which $A$ and $B$ are defined. We shall also use the usual bracket $[A,\,B]:=AB -(-1)^{ab}\,BA$ for graded linear operators $A, B$ of respective degrees $a, b$ on the algebra $\oplus_k\Lambda^kT^\star X$ of differential forms on $X$.}

\end{Not}

\noindent {\bf Acknowledgments.} The author wishes to thank L. Ugarte for useful discussions about the content of this paper and for suggestions for section \ref{section:consequences_MT1}. Thanks are also due to S. Rao and Q. Zhao for stimulating discussions.

\section{Rescaled Laplacians}\label{section:rescaled-laplacians} Let $X$ be a compact complex manifold with $\mbox{dim}_\C X=n$. We fix a Hermitian metric $\omega$ on $X$.

\subsection{Rescaling the metric}\label{subsection:rescaling_metric} The first operation we will consider is a {\it partial rescaling of $\omega$} in a way that depends solely on the {\it holomorphic} degree of forms.

\begin{Def}\label{Def:metric-rescaling} For all $p,q\in\{0,\dots , n\}$, all $(p,\,q)$-forms $u,v$ and every constant $h>0$, we define the following {\bf pointwise inner product}

$$\langle u,\,v\rangle_{\omega_h}:=h^{2p}\,\langle u,\,v\rangle_\omega$$

\noindent where $\langle,\,\,\,\rangle_\omega$ stands for the pointwise inner product defined by the original Hermitian metric $\omega$.

\end{Def}

Note that, for every $h>0$, we obtain in this way a Hermitian metric $\omega_h$ on every vector bundle $\Lambda^{p,\,q}T^\star X$ of $(p,\,q)$-forms on $X$. The maps

$$\theta_h:\Lambda^{p,\,q}T^\star X\longrightarrow\Lambda^{p,\,q}T^\star X, \hspace{3ex} u\mapsto\theta_h u:=h^pu,$$

\noindent induce an {\bf isometry} of Hermitian vector bundles $\theta_h:(\Lambda T^\star X,\,\omega_h)\longrightarrow(\Lambda T^\star X,\,\omega)$ since 

$$\langle u,\,v\rangle_{\omega_h} = \langle h^pu,\,h^pv\rangle_\omega = \langle \theta_hu,\,\theta_hv\rangle_\omega \hspace{3ex} \mbox{for all}\hspace{1ex} u,v\in\Lambda^{p,\,q}T^\star X.$$

In particular, we have defined a Hermitian metric 

$$\omega_h=\frac{1}{h^2}\,\omega, \hspace{3ex} h>0,$$ 

\noindent on the holomorphic tangent bundle $T^{1,\,0}X$ of vector fields of type $(1,\,0)$, or equivalently, a rescaled $C^\infty$ positive-definite $(1,\,1)$-form $\omega_h=h^{-2}\,\omega$ on $X$. This induces a $C^\infty$ positive volume form 

$$dV_{\omega_h}:=\frac{\omega_h^n}{n!} = \frac{1}{h^{2n}}\,\frac{\omega^n}{n!}=\frac{1}{h^{2n}}\,dV_\omega$$

\noindent on $X$, which in turn gives rise, in conjunction with the above pointwise inner product $\langle\,\,,\,\,\rangle_{\omega_h}$, to the following {\bf $L^2$ inner product}

$$\langle\langle u,\, v\rangle\rangle_{\omega_h}:=\int\limits_X \langle u,\,v\rangle_{\omega_h}\,dV_{\omega_h} = \frac{1}{h^{2n}}\,\int\limits_X \langle\theta_h u,\,\theta_h v\rangle_\omega\,dV_\omega = \frac{1}{h^{2n}}\,\langle\langle\theta_h u,\,\theta_h v\rangle\rangle_\omega$$

\noindent for all forms $u,v\in C^\infty_{p,\,q}(X,\,\C)$ and all bidegrees $(p,\,q)$.

\begin{Formula}\label{Formula:L2_comparison} For all $(p,\,q)$-forms $u,v$, we have

$$\langle\langle u,\, v\rangle\rangle_{\omega_h} = \frac{1}{h^{2(n-p)}}\,\langle\langle u,\, v\rangle\rangle_\omega, \hspace{3ex} \mbox{hence} \hspace{3ex} ||u||_{\omega_h} = h^{-(n-p)}\,||u||_\omega.$$

\end{Formula}

\noindent {\it Proof.} The formula follows at once from the last identity and from the fact that $\theta_h u = h^pu$ for all $(p,\,q)$-forms $u$. \hfill $\Box$

\begin{Def}\label{Def:omega_h_Laplacian} Let $(X,\,\omega)$ be a compact Hermitian manifold with $\mbox{dim}_\C X=n$. For every $k=0,\dots , 2n$ and every constant $h>0$, we consider the {\bf $d$-Laplacian w.r.t. the rescaled metric $\omega_h$} acting on $C^\infty$ $k$-forms on $X$:

$$\Delta_{\omega_h}:C^\infty_k(X,\,\C)\longrightarrow C^\infty_k(X,\,\C), \hspace{3ex} \Delta_{\omega_h}:=dd^\star_{\omega_h} + d^\star_{\omega_h}d,$$

\noindent where $d^\star_{\omega_h}$ is the formal adjoint of $d$ w.r.t. $\langle\langle\,\,,\,\,\rangle\rangle_{\omega_h}$ and $\langle\langle\,\,,\,\,\rangle\rangle_{\omega_h}$ has been extended from the spaces $C^\infty_{p,\,q}(X,\,\C)$ to $C^\infty_k(X,\,\C)=\oplus_{p+q=k}C^\infty_{p,\,q}(X,\,\C)$ by sesquilinearity and by imposing that $\langle\langle u,\,v\rangle\rangle_{\omega_h}=0$ whenever $u\in C^\infty_{p,\,q}(X,\,\C)$ and $v\in C^\infty_{r,\,s}(X,\,\C)$ with $(p,\,q)\neq(r,\,s)$.

\end{Def}

\subsection{Rescaling the differential}\label{subsection:rescaling_metric} The second operation we will consider is a {\it partial rescaling of $d=\partial + \bar\partial$} that applies solely to its component of type $(1,\,0)$.

\begin{Def}\label{Def:d-rescaling} Let $X$ be a compact complex manifold, $\mbox{dim}_\C X=n$. For every constant $h>0$, let 

$$d_h:=h\partial + \bar\partial:C^\infty_k(X,\,\C)\longrightarrow C^\infty_{k+1}(X,\,\C), \hspace{3ex} k\in\{0,\,\dots , 2n\}.$$

\end{Def}

Some basic properties of the rescaled differential $d_h$ are summed up in the following

\begin{Lem}\label{Lem:d_h_properties} $(i)$ The operators $d$ and $d_h$ are related by the identity

$$d_h =\theta_hd\theta_h^{-1}.$$

\noindent $(ii)$\, $d_h^2=0$ and the $d$- and $d_h$-cohomologies are related by the {\bf isomorphism}

$$H^k_d(X,\,\C)\stackrel{\simeq}{\longrightarrow} H^k_{d_h}(X,\,\C), \hspace{3ex} \{u\}_d\mapsto \{\theta_hu\}_{d_h}$$

\noindent where $H^k_d(X,\,\C) = H^k_{DR}(X,\,\C)$ are the usual De Rham cohomology groups, while $H^k_{d_h}(X,\,\C):=\ker(d_h:C^\infty_k(X,\,\C)\longrightarrow C^\infty_{k+1}(X,\,\C))/\mbox{Im}\,(d_h:C^\infty_{k-1}(X,\,\C)\longrightarrow C^\infty_k(X,\,\C))$ are the $d_h$-cohomology groups.

\end{Lem}

\noindent {\it Proof.} $(i)$\, If $u$ is a $(p,\,q)$-form, we have

\vspace{1ex}

$(\theta_hd\theta_h^{-1})(u) = \theta_hd(h^{-p}u) = h^{-p}\theta_h(\partial u) + h^{-p}\theta_h(\bar\partial u) = h^{-p}h^{p+1}\partial u + h^{-p}h^p\bar\partial u = h\partial u + \bar\partial u = d_hu.$

\vspace{1ex}

\noindent Thus, $d_h =\theta_hd\theta_h^{-1}$ on pure-type forms, so this identity extends to arbitrary forms by linearity.

$(ii)$\, On the one hand, $d_h^2=\theta_hd^2\theta_h^{-1}=0$; on the other hand, 

\vspace{1ex}

\hspace{1ex} $d_h(\theta_hu) = \theta_hdu$, \hspace{3ex} so we have the equivalence: \hspace{3ex} $\theta_hu\in\ker(d_h) \iff u\in\ker d$;

\vspace{1ex}

\hspace{1ex} $\theta_hu = d_hv$ iff $u=d(\theta_h^{-1}v)$, \hspace{3ex} so we have the equivalence: \hspace{3ex} $\theta_hu\in\mbox{Im}\,(d_h) \iff u\in\mbox{Im}\, d$.

\vspace{2ex}

\noindent These equivalences show that the linear map $H^k_d(X,\,\C)\ni \{u\}_d\mapsto \{\theta_hu\}_{d_h}\in H^k_{d_h}(X,\,\C)$ is well defined and bijective. \hfill $\Box$

\vspace{2ex}

 In particular, the spectral sequences induced by the pairs of differentials $(\partial,\,\bar\partial)$ and $(h\partial,\,\bar\partial)$ are {\it isomorphic}, so degenerate at the same page. The first of them is the Fr\"olicher spectral sequence of $X$.

\begin{Def}\label{Def:h_Laplacian} Let $(X,\,\omega)$ be a compact Hermitian manifold with $\mbox{dim}_\C X=n$. For every constant $h>0$ and every degree $k\in\{0,\dots, 2n\}$, we consider the {\bf $d_h$-Laplacian w.r.t. the given metric $\omega$} acting on $C^\infty$ $k$-forms on $X$:

 $$\Delta_h:C^\infty_k(X,\,\C)\longrightarrow C^\infty_k(X,\,\C), \hspace{3ex} \Delta_h:=d_hd_h^\star + d_h^\star d_h,$$

\noindent where $d_h^\star$ is the formal adjoint of $d_h$ w.r.t. the $L^2$ inner product induced by $\omega$.
 
\end{Def}

\subsection{Comparison of the two rescaled Laplacians}\label{subsection:rescaled_Laplacians}

We now bring together the above two operations by comparing the corresponding Laplace-type operators. Note that $\Delta_{\omega_h}$ was defined by the rescaled differential $d_h$ and the original metric $\omega$, while $\Delta_h$ was induced by the rescaled metric $\omega_h$ and the original differential $d$. 

\begin{Lem}\label{Lem:relations_Laplacians} $(i)$\, If $\theta_h^\star$ and $d_h^\star$ stand for the formal adjoints of $\theta_h$, resp. $d_h$, w.r.t. the pointwise, resp. $L^2$, inner product induced by $\omega$, we have

$$\theta_h^\star = \theta_h \hspace{3ex} \mbox{and} \hspace{3ex} d_h^\star = \theta_h^{-1}d^\star\theta_h.$$

\noindent $(ii)$\, The adjoints $\partial^\star_{\omega_h}$, $\bar\partial^\star_{\omega_h}$ w.r.t. to the metric $\omega_h$, as well as the adjoints $\partial^\star_\omega=\partial^\star$, $\bar\partial^\star_\omega=\bar\partial^\star$ w.r.t. to the metric $\omega$, of $\partial$, resp. $\bar\partial$ are related by the formulae:

$$\partial^\star_{\omega_h} = h^2\partial^\star  \hspace{3ex} \mbox{and} \hspace{3ex} \bar\partial^\star_{\omega_h} = \bar\partial^\star.$$

\noindent Consequently, we get \begin{eqnarray}\nonumber\Delta_{\omega_h} & = & h^2\Delta' + \Delta'' + [\partial,\,\bar\partial^\star] + h^2[\bar\partial,\,\partial^\star] \\
\nonumber & = & h^2\Delta' + \Delta'' - [\partial,\,\bar\tau^\star] - h^2[\bar\tau,\,\partial^\star] = h^2\Delta' + \Delta'' - [\tau,\,\bar\partial^\star] - h^2[\bar\partial,\,\tau^\star],\end{eqnarray} \noindent and \begin{eqnarray}\nonumber\Delta_h & = & h^2\Delta' + \Delta'' + h[\partial,\,\bar\partial^\star] + h[\bar\partial,\,\partial^\star] \\
\nonumber & = & h^2\Delta' + \Delta'' - h[\partial,\,\bar\tau^\star] - h[\bar\tau,\,\partial^\star] = h^2\Delta' + \Delta'' - h[\tau,\,\bar\partial^\star] - h[\bar\partial,\,\tau^\star],\end{eqnarray}

\noindent where the adjoints $\partial^\star, \bar\partial^\star, \tau^\star, \bar\tau^\star$ and the Laplacians $\Delta', \Delta''$ are computed w.r.t. the metric $\omega$, while

$$\tau=\tau_\omega:=[\Lambda_\omega,\,\partial\omega\wedge\cdot]: C^\infty_{p,\,q}(X,\,\C)\longrightarrow C^\infty_{p+1,\,q}(X,\,\C)$$

\noindent is the {\bf torsion operator} (of type $(1,\,0)$ and order zero, acting on smooth forms of any bidegree $(p,\,q)$, where $\Lambda_\omega$ is the adjoint of the multiplication operator $\omega\wedge\cdot$) associated with the metric $\omega$ as defined in $[Dem84]$ (see also $[Dem97, VII, \S.1]$).

 In particular, the second-order Laplacians $\Delta_{\omega_h}$ and $\Delta_h$ are {\bf elliptic} since the second-order Laplacians $\Delta'$ and $\Delta''$ are and the deviation terms $- [\partial,\,\bar\tau^\star] - h^2[\bar\tau,\,\partial^\star]$ and $-h[\partial,\,\bar\tau^\star] - h[\bar\tau,\,\partial^\star] $ are only of order $1$. 

Note that $\langle\langle[\partial,\,\bar\partial^\star]u,\,u\rangle\rangle = \langle\langle[\bar\partial,\,\partial^\star]u,\,u\rangle\rangle = 0$ whenever the form $u$ is {\bf of pure type} and whatever metric is used to define $\langle\langle\,\,,\,\,\rangle\rangle$ (because pure-type forms of different bidegrees are orthogonal w.r.t. any metric), so \begin{equation}\label{eqn:pure-type_vanishing}\langle\langle\Delta_{\omega_h}u,\,u\rangle\rangle = \langle\langle\Delta_h u,\,u\rangle\rangle = h^2\,\langle\langle\Delta' u,\,u\rangle\rangle + \langle\langle\Delta'' u,\,u\rangle\rangle \hspace{3ex} \mbox{for every} \hspace{1ex} \mbox{\bf pure-type} \hspace{1ex} \mbox{form}\hspace{1ex} u.\end{equation}

(This fails, in general, if $u$ is not of pure type, unless the metric $\omega$ is K\"ahler.)

\vspace{1ex}

\noindent $(iii)$\, The rescaled Laplacians $\Delta_{\omega_h}$ and $\Delta_h$ are related by the formula \begin{eqnarray}\label{eqn:rescaled-Laplacians_relation}\Delta_h = \theta_h\Delta_{\omega_h}\theta_h^{-1}.\end{eqnarray}

\end{Lem}

\noindent {\it Proof.} $(i)$\, For any $k$-forms $u=\sum\limits_{p+q=k}u^{p,\,q}$ and $v=\sum\limits_{p+q=k}v^{p,\,q}$, we have 

\vspace{1ex}

$\langle\theta_hu,\,v\rangle_\omega = \sum\limits_{p+q=k}\langle h^pu^{p,\,q},\,v^{p,\,q}\rangle_\omega = \sum\limits_{p+q=k}\langle u^{p,\,q},\,h^pv^{p,\,q}\rangle_\omega = \langle u,\,\theta_hv\rangle_\omega,$ \hspace{2ex} so $\theta_h^\star = \theta_h$.

\vspace{1ex}

\noindent The second identity in $(i)$ follows by taking conjugates in $d_h =\theta_hd\theta_h^{-1}$.

$(ii)$\, For any forms $\alpha\in C^\infty_{p-1,\,q}(X,\,\C)$ and $\beta\in C^\infty_{p,\,q}(X,\,\C)$, we have

\begin{eqnarray}\nonumber\langle\langle\alpha,\,\partial^\star_\omega\beta\rangle\rangle_\omega & = & \langle\langle\partial\alpha,\,\beta\rangle\rangle_\omega = \int\limits_X\langle\partial\alpha,\,\beta\rangle_\omega\,dV_\omega = \int\limits_X\frac{1}{h^{2p}}\,\langle\partial\alpha,\,\beta\rangle_{\omega_h}\,h^{2n}dV_{\omega_h} = h^{2(n-p)}\,\langle\langle\partial\alpha,\,\beta\rangle\rangle_{\omega_h}\\
\nonumber & = & h^{2(n-p)}\,\langle\langle\alpha,\,\partial^\star_{\omega_h}\beta\rangle\rangle_{\omega_h} = h^{2(n-p)}\,\int\limits_Xh^{2(p-1)}\,\langle\alpha,\,\partial^\star_{\omega_h}\beta\rangle_\omega\,\frac{1}{h^{2n}}\,dV_\omega = \frac{1}{h^2}\,\langle\langle\alpha,\,\partial^\star_{\omega_h}\beta\rangle\rangle_\omega.\end{eqnarray}

\noindent We get $\partial^\star_\omega=h^{-2}\,\partial^\star_{\omega_h}$, which is the first identity under $(ii)$. 

 The identity $\bar\partial^\star_{\omega_h} = \bar\partial^\star_\omega$ is proved in the same way by using the fact that $\bar\partial$ acts only on the anti-holomorphic degree of forms which is unaffected by the change of metric from $\omega$ to $\omega_h$. 

 Using these formulae, we get \begin{eqnarray}\nonumber\Delta_{\omega_h} & = & [\partial+\bar\partial,\,\partial^\star_{\omega_h} + \bar\partial^\star_{\omega_h}] = [\partial,\, h^2\partial^\star] + [\bar\partial,\, \bar\partial^\star] + [\partial,\, \bar\partial^\star] + [\bar\partial,\, h^2\partial^\star] \\
\nonumber & = & h^2\Delta' + \Delta'' +  [\partial,\, \bar\partial^\star] + h^2[\bar\partial,\, \partial^\star]\end{eqnarray}  and \begin{eqnarray}\nonumber\Delta_h & = & [h\partial+\bar\partial,\,h\partial^\star + \bar\partial^\star] = h^2[\partial,\,\partial^\star] + [\bar\partial,\, \bar\partial^\star] + h[\partial,\, \bar\partial^\star] + h[\bar\partial,\, \partial^\star] \\
\nonumber & = & h^2\Delta' + \Delta'' +  h[\partial,\, \bar\partial^\star] + h[\bar\partial,\, \partial^\star].\end{eqnarray}

 On the other hand, we know from [Dem84] (or [Dem97, VII, $\S.1$]) that $$[\partial,\,\bar\partial^\star] = -[\partial,\,\bar\tau^\star] = -[\tau,\,\bar\partial^\star] \hspace{3ex} \mbox{and, by conjugation, we get} \hspace{3ex} [\bar\partial,\,\partial^\star] = -[\bar\partial,\,\tau^\star] = -[\bar\tau,\,\partial^\star].$$

\noindent So, the terms measuring the deviations of $\Delta_{\omega_h}$ and $\Delta_h$ from $h^2\Delta' + \Delta''$ are of order $1$ and we get the alternative formulae for $\Delta_{\omega_h}$ and $\Delta_h$ spelt out in the statement.

$(iii)$\, For any smooth $(p,\,q)$-form $\alpha$, we have \begin{eqnarray}\nonumber(\theta_h\Delta_{\omega_h}\theta_h^{-1})\alpha & = & \frac{1}{h^p}\theta_h\Delta_{\omega_h}\alpha = \frac{1}{h^p}\theta_h(h^2\Delta'\alpha) + \frac{1}{h^p}\theta_h(\Delta''\alpha) + \frac{1}{h^p}\theta_h([\partial,\,\bar\partial^\star]\alpha) + \frac{1}{h^p}\theta_h(h^2[\bar\partial,\,\partial^\star]\alpha) \\
\nonumber & = &  \frac{h^2h^p}{h^p}\Delta'\alpha + \frac{h^p}{h^p}\Delta''\alpha + \frac{h^{p+1}}{h^p}[\partial,\,\bar\partial^\star]\alpha + \frac{h^2h^{p-1}}{h^p}[\bar\partial,\,\partial^\star]\alpha \\
\nonumber & = & h^2\Delta'\alpha + \Delta''\alpha + h[\partial,\,\bar\partial^\star]\alpha + h[\bar\partial,\,\partial^\star]\alpha = \Delta_h\alpha.\end{eqnarray}

\noindent Thus, $\theta_h\Delta_{\omega_h}\theta_h^{-1} = \Delta_h$ on pure-type forms and this identity extends to arbitrary forms by linearity.   \hfill $\Box$

\begin{Cor}\label{Cor:spectra} Let $(X,\,\omega)$ be a compact Hermitian manifold with $\mbox{dim}_\C X=n$. For every constant $h>0$ and every degree $k\in\{0,\dots, 2n\}$, the {\bf spectra} of the rescaled Laplacians $\Delta_h, \Delta_{\omega_h}:C^\infty_k(X,\,\C)\longrightarrow C^\infty_k(X,\,\C)$ {\bf coincide}, i.e.

\begin{equation}\label{eqn:spectra_equality}\mbox{Spec}(\Delta_h) = \mbox{Spec}(\Delta_{\omega_h}),\end{equation}

\noindent and their respective {\bf eigenspaces} are obtained from each other via the rescaling isometry $\theta_h$:

\begin{equation}\label{eqn:eigenspaces_transformation}\theta_h(E_{\Delta_{\omega_h}}(\lambda)) = E_{\Delta_h}(\lambda) \hspace{6ex} \mbox{for every} \hspace{1ex} \lambda\in\mbox{Spec}(\Delta_h) = \mbox{Spec}(\Delta_{\omega_h}),    \end{equation}

\noindent where $E_{\Delta_{\omega_h}}(\lambda)$, resp. $E_{\Delta_h}(\lambda)$, stands for the eigenspace corresponding to the eigenvalue $\lambda$ of the operator $\Delta_{\omega_h}$, resp. $\Delta_h$. 

Thus, $\Delta_h$ and $\Delta_{\omega_h}$ have the {\bf same eigenvalues} with the {\bf same multiplicities}.

\end{Cor}

\noindent {\it Proof.} Let $\lambda\in\mbox{Spec}(\Delta_{\omega_h})$ and let $\alpha\in E_{\Delta_{\omega_h}}(\lambda)\subset C^\infty_k(X,\,\C)$. So $\Delta_{\omega_h}\alpha = \lambda\alpha$, hence

$$\Delta_h(\theta_h\alpha) = (\theta_h\Delta_{\omega_h}\theta_h^{-1})(\theta_h\alpha) = \theta_h(\lambda\alpha) = \lambda(\theta_h\alpha).$$

\noindent Thus, $\lambda\in\mbox{Spec}(\Delta_h)$ and $\theta_h\alpha\in E_{\Delta_h}(\lambda)$. These implications also hold in reverse order, so we get the equivalences:

$$\lambda\in\mbox{Spec}(\Delta_h) \iff \lambda\in\mbox{Spec}(\Delta_{\omega_h}) \hspace{3ex} \mbox{and} \hspace{3ex} \alpha\in E_{\Delta_{\omega_h}}(\lambda) \iff \theta_h\alpha\in E_{\Delta_h}(\lambda).$$

 These equivalences amount to (\ref{eqn:spectra_equality}) and (\ref{eqn:eigenspaces_transformation}).   \hfill $\Box$

\vspace{3ex}

Another consequence of the above discussion is a Hodge Theory for the $d_h$-cohomology and the resulting equidimensionality of the kernels of $\Delta$ and $\Delta_h$ in every degree.

\begin{Cor}\label{Cor:Hodge_d_h} Let $(X,\,\omega)$ be a compact Hermitian manifold with $\mbox{dim}_\C X=n$. For every constant $h>0$ and every degree $k\in\{0,\dots, 2n\}$, the operator $d_h:C^\infty_k(X,\,\C)\longrightarrow C^\infty_k(X,\,\C)$ induces the following $L^2_\omega$-orthogonal direct-sum decomposition:

$$C^\infty_k(X,\,\C) = {\cal H}^k_{\Delta_h}(X,\,\C)\oplus\mbox{Im}\,d_h\oplus\mbox{Im}\,d_h^\star,$$
 
\noindent where ${\cal H}^k_{\Delta_h}(X,\,\C)$ is the kernel of $\Delta_h:C^\infty_k(X,\,\C)\longrightarrow C^\infty_k(X,\,\C)$ and $\ker d_h={\cal H}^k_{\Delta_h}(X,\,\C)\oplus\mbox{Im}\,d_h$. The vector space ${\cal H}^k_{\Delta_h}(X,\,\C)$ is {\bf finite-dimensional}, while $\mbox{Im}\,d_h$ and $\mbox{Im}\,d_h^\star$ are {\bf closed} subspaces of $C^\infty_k(X,\,\C)$.

This, in turn, induces the {\bf Hodge isomorphism}

$${\cal H}^k_{\Delta_h}(X,\,\C) \simeq H^k_{d_h}(X,\,\C), \hspace{3ex} \alpha\mapsto\{\alpha\}_{d_h}.$$

 Since $H^k_d(X,\,\C)$ and $H^k_{d_h}(X,\,\C)$ are isomorphic (via $\theta_h$, see Lemma \ref{Lem:d_h_properties}) and ${\cal H}^k_\Delta(X,\,\C)\simeq H^k_d(X,\,\C)$ (by standard Hodge theory), we infer that ${\cal H}^k_\Delta(X,\,\C)$ and ${\cal H}^k_{\Delta_h}(X,\,\C)$ are {\bf isomorphic} (although the isomorphism need not be defined by $\theta_h$). In particular,

$$\mbox{dim}\,{\cal H}^k_{\Delta_h}(X,\,\C) = \mbox{dim}\,{\cal H}^k_\Delta(X,\,\C) \hspace{3ex} \mbox{for all} \hspace{1ex} h>0.$$

\end{Cor}

\noindent {\it Proof.} Since $X$ is compact and $\Delta_h$ is elliptic and self-adjoint, a standard consequence of G${\mathring a}$rding's inequality (see e.g. [Dem97, VI]) yields the two-space orthogonal decomposition $C^\infty_k(X,\,\C) = {\cal H}^k_{\Delta_h}(X,\,\C)\oplus\mbox{Im}\,\Delta_h$, while this, together with the integrability property $d_h^2=0$, further induces the orthogonal splitting $\mbox{Im}\,\Delta_h = \mbox{Im}\,d_h\oplus\mbox{Im}\,d_h^\star$. The same consequence of G${\mathring a}$rding's inequality ensures that $\ker\Delta_h$ is finite-dimensional and that the images in $C^\infty_k(X,\,\C)$ of $d_h$ and $d_h^\star$ are closed. \hfill $\Box$

\section{The differentials in the Fr\"olicher spectral sequence}\label{section:differentials-Froelicher}

We begin by recalling the well-known construction of the Fr\"olicher spectral sequence in order to fix the notation and to point out the key features for us.

Let $X$ be a compact complex manifold with $\mbox{dim}_\C X=n$. Recall that the zero-th page $E_0$ of the Fr\"olicher spectral sequence of $X$ consists of the spaces $E_0^{p,\,q}:=C^\infty_{p,\,q}(X,\,\C)$ of smooth pure-type forms on $X$ and of the type-$(0,\,1)$ differentials $d_0:=\bar\partial$ forming the Dolbeault complex:

$$\dots\stackrel{d_0}{\longrightarrow}E_0^{p,\,q-1}\stackrel{d_0}{\longrightarrow}E_0^{p,\,q}\stackrel{d_0}{\longrightarrow}E_0^{p,\,q+1}\stackrel{d_0}{\longrightarrow}\dots.$$ 

\noindent Thus, in every bidegree $(p,\,q)$, the inclusions $\mbox{Im}\,d_0^{p,\,q-1}\subset\ker d_0^{p,\,q}\subset E_0^{p,\,q}$ induce (infinitely many, non-canonical) isomorphisms

\begin{equation}\label{eqn:E_0_splitting}C^\infty_{p,\,q}(X,\,\C) \simeq \mbox{Im}\,d_0^{p,\,q-1}\oplus E_1^{p,\,q}\oplus(E_0^{p,\,q}/\ker d_0^{p,\,q}),\end{equation}

\noindent where $d_0=d_0^{p,\,q}:E_0^{p,\,q}\longrightarrow E_0^{p,\,q+1}$ is the differential $d_0$ acting in bidegree $(p,\,q)$ and the $E_1^{p,\,q}:=\ker d_0^{p,\,q}/\mbox{Im}\,d_0^{p,\,q-1}=H^{p,\,q}_{\bar\partial}(X,\,\C)$ are the Dolbeault cohomology groups of $X$.

The first page $E_1$ of the Fr\"olicher spectral sequence consists of the spaces $E_1^{p,\,q}$ (i.e. the cohomology of the zero-th page) and of the type-$(1,\,0)$ differentials $d_1$:

$$\dots\stackrel{d_1}{\longrightarrow}E_1^{p-1,\,q}\stackrel{d_1}{\longrightarrow}E_1^{p,\,q}\stackrel{d_1}{\longrightarrow}E_1^{p+1,\,q}\stackrel{d_1}{\longrightarrow}\dots.$$

\noindent induced in cohomology by $\partial$ (i.e. $d_1([\alpha]_{\bar\partial}):=[\partial\alpha]_{\bar\partial}$). Thus, in every bidegree $(p,\,q)$, the inclusions $\mbox{Im}\,d_1^{p-1,\,q}\subset\ker d_1^{p,\,q}\subset E_1^{p,\,q}$ induce (infinitely many, non-canonical) isomorphisms

\begin{equation}\label{eqn:E_1_splitting}E_1^{p,\,q} \simeq \mbox{Im}\,d_1^{p-1,\,q}\oplus E_2^{p,\,q}\oplus(E_1^{p,\,q}/\ker d_1^{p,\,q}),\end{equation}

\noindent where $d_1^{p,\,q}$ is $d_1$ acting in bidegree $(p,\,q)$, while the spaces $E_2^{p,\,q}:=\ker d_1^{p,\,q}/\mbox{Im}\,d_1^{p-1,\,q}$ form the cohomology of the page $E_1$.

 The remaining pages are constructed inductively: the differentials $d_r=d_r^{p,\,q}:E_r^{p,\,q}\longrightarrow E_r^{p+r,\,q-r+1}$ are of type $(r,\,-r+1)$ for every $r$, while the spaces $E_r^{p,\,q}:=\ker d_{r-1}^{p,\,q}/\mbox{Im}\,d_{r-1}^{p-r+1,\,q+r-2}$ on the $r^{th}$ page are defined as the cohomology of the previous page $E_{r-1}$. On every page $E_r$ and in every bidegree $(p,\,q)$, the inclusions $\mbox{Im}\,d_r^{p-r,\,q+r-1}\subset\ker d_r^{p,\,q}\subset E_r^{p,\,q}$ induce (infinitely many, non-canonical) isomorphisms

\begin{equation}\label{eqn:E_r_splitting}E_r^{p,\,q} \simeq \mbox{Im}\,d_r^{p-r,\,q+r-1}\oplus E_{r+1}^{p,\,q}\oplus(E_r^{p,\,q}/\ker d_r^{p,\,q}),\end{equation}

\noindent where $E_{r+1}^{p,\,q}:=\ker d_r^{p,\,q}/\mbox{Im}\,d_r^{p-r,\,q+r-1}$.

It is worth stressing that (\ref{eqn:E_0_splitting}), (\ref{eqn:E_1_splitting}) and (\ref{eqn:E_r_splitting}) only assert that the vector spaces on either side of $\simeq$ are isomorphic, but no choice of preferred isomorphism is possible at this stage.

 A classical result of Fr\"olicher [Fro55] asserts that this spectral sequence converges to the De Rham cohomology of $X$ and degenerates after finitely many steps. This means that there are (non-canonical) isomorphisms:

\begin{equation}\label{eqn:Froelicher_decomposition}H^k_{DR}(X,\,\C)\simeq\bigoplus\limits_{p+q=k}E_\infty^{p,\,q}, \hspace{6ex} k=0,\dots , 2n,\end{equation}

\noindent where $E_\infty^{p,\,q} = \dots = E_{N+2}^{p,\,q} = E_{N+1}^{p,\,q}= E_N^{p,\,q}$ for all $p,q$ and where $N\geq 1$ is the positive integer such that the spectral sequence degenerates at $E_N$.

\subsection{Identification of the $d_r$'s with restrictions of $d$}\label{subsection:d_r-d_identification}

 Summing up (\ref{eqn:E_0_splitting}), (\ref{eqn:E_1_splitting}), (\ref{eqn:E_r_splitting}) over $r=0,\dots , N-1$, we get (infinitely many, non-canonical) isomorphisms

$$C^\infty_{p,\,q}(X,\,\C) \simeq \bigoplus\limits_{r=0}^{N-1}\mbox{Im}\,d_r^{p-r,\,q+r-1} \oplus E_\infty^{p,\,q} \oplus \bigoplus\limits_{r=0}^{N-1}(E_r^{p,\,q}/\ker d_r^{p,\,q}) $$

\noindent for every bidegree $(p,\,q)$. Note that the isomorphisms (\ref{eqn:E_0_splitting}), (\ref{eqn:E_1_splitting}), (\ref{eqn:E_r_splitting}) identify the spaces $\mbox{Im}\,d_r^{p-r,\,q+r-1}$, $E_r^{p,\,q}$ (including for $r=\infty$) and $E_r^{p,\,q}/\ker d_r^{p,\,q}$ with certain subspaces of $C^\infty_{p,\,q}(X,\,\C)$. However, these subspaces have not been specified yet since multiple choices (and no canonical choice) are possible for the isomorphisms (\ref{eqn:E_0_splitting}), (\ref{eqn:E_1_splitting}), (\ref{eqn:E_r_splitting}). These choices can only be made unique once a Hermitian metric has been fixed on $X$. (See $\S.$\ref{subsection:explicit}.)

Now, since $C^\infty_k(X,\,\C)=\oplus_{p+q=k}C^\infty_{p,\,q}(X,\,\C)$ for all $k=0,\dots , 2n$, we get

$$\xymatrix{
 C^\infty_k(X,\,\C)\ar[d]^d \simeq & \bigoplus\limits_{0\leq r\leq N-1\atop p+q=k}\mbox{Im}\,d_r^{p-r,\,q+r-1} & \hspace{-20ex} \oplus \bigoplus\limits_{p+q=k}E_\infty^{p,\,q} \oplus \bigoplus\limits_{0\leq r\leq N-1\atop p+q=k}(E_r^{p,\,q}/\ker d_r^{p,\,q})  \\
 C^\infty_{k+1}(X,\,\C) \simeq & \bigoplus\limits_{0\leq r\leq N-1\atop p+q=k}\mbox{Im}\,d_r^{p,\,q} & \hspace{-8ex} \oplus \bigoplus\limits_{p'+q'=k+1}E_\infty^{p',\,q'} \oplus  \bigoplus\limits_{0\leq r\leq N-1\atop p+q=k}(E_r^{p+r,\,q-r+1}/\ker d_r^{p+r,\,q-r+1}).   }
$$

\noindent Thus, under these isomorphisms, the operator $d=d^{(k)}:C^\infty_k(X,\,\C)\longrightarrow C^\infty_{k+1}(X,\,\C)$ identifies as

\begin{equation}\label{eqn:d_splitting}d^{(k)}\simeq\bigoplus\limits_{0\leq r\leq N-1\atop p+q=k}d_r^{p,\,q},\end{equation} 

\noindent where the isomorphism $d_r^{p,\,q}:E_r^{p,\,q}/\ker d_r^{p,\,q}\longrightarrow \mbox{Im}\,d_r^{p,\,q}$ is the restriction of $d_r=d_r^{p,\,q}:E_r^{p,\,q}\longrightarrow E_r^{p+r,\,q-r+1}$ to the third piece on the r.h.s. of (\ref{eqn:E_r_splitting}). The fact that $d_r$ is of type $(r,\,-r+1)$ will play a key role in the sequel. 

 On the other hand, summing up the splittings of $C^\infty_{p,\,q}(X,\,\C)$ over $p\geq s$ for any given $s$, we get

$${\cal A}_s^k:=\bigoplus\limits_{p\geq s\atop p+q=k} C^\infty_{p,\,q}(X,\,\C) \simeq \bigoplus\limits_{p\geq s\atop p+q=k}\bigg[\bigoplus\limits_{r=0}^{N-1}\mbox{Im}\,d_r^{p-r,\,q+r-1} \oplus E_\infty^{p,\,q} \oplus \bigoplus\limits_{r=0}^{N-1}(E_r^{p,\,q}/\ker d_r^{p,\,q})\bigg].$$

\begin{Lem}\label{Lem:E_r^k_dimension} $(i)$\, For every $r$ and every $k$, let $E_r^k:=\bigoplus\limits_{p+q=k}E_r^{p,\,q}$. Then

\begin{equation}\label{eqn:dim_E_r^k}\mbox{dim}\,E_r^k = \sum\limits_{p+q=k}\mbox{dim}\,E_r^{p,\,q} = b_k + m_r^{k-1} + m_r^k, \hspace{6ex} 0\leq r\leq N, \hspace{1ex} 0\leq k\leq 2n,\end{equation}

\noindent where we set $m_r^k:=\sum\limits_{l\geq r}\sum\limits_{p+q=k}\mbox{dim}\,(E_l^{p,\,q}/\ker d_l^{p,\,q})$.

 $(ii)$\, For every $r$ and every $k$, let $L_r^{p,\,q}:=\bigoplus\limits_{l\geq r}(E_l^{p,\,q}/\ker d_l^{p,\,q})$ and $L_r^k:=\bigoplus\limits_{p+q=k}L_r^{p,\,q}$. Then, $\mbox{dim}\,L_r^k=m_r^k$ (obvious) and, under the identifications defined by the isomorphisms (\ref{eqn:E_0_splitting}), (\ref{eqn:E_1_splitting}), (\ref{eqn:E_r_splitting}), the following inclusions hold:  

\begin{equation}\label{eqn:d_L_r_pq}d(L_r^{p,\,q})\subset{\cal A}_{p+r}^{p+q+1}, \hspace{6ex} 0\leq r\leq N, \hspace{1ex} 0\leq p,q\leq n,\end{equation}

\noindent where $d(L_r^{p,\,q}):=\oplus_{l\geq r}d_l^{p,\,q}(E_l^{p,\,q}/\ker d_l^{p,\,q})$ in keeping with identification (\ref{eqn:d_splitting}).

\end{Lem}

\noindent {\it Proof.} $(i)$\, For every fixed $r$, summing up the splittings (\ref{eqn:E_r_splitting}) with $l$ in place of $r$ over $l\geq r$ and then summing up over $p+q=k$ for every fixed $k$, we get

\begin{eqnarray}\nonumber E_r^k\simeq \bigoplus\limits_{p+q=k}E_\infty^{p,\,q} \oplus \bigoplus\limits_{l\geq r}\bigoplus\limits_{p+q=k}\mbox{Im}\,d_l^{p-l,\,q+l-1} \oplus\bigoplus\limits_{l\geq r}\bigoplus\limits_{p+q=k}(E_l^{p,\,q}/\ker d_l^{p,\,q}).\end{eqnarray}

\noindent Since $\mbox{Im}\,d_l^{p-l,\,q+l-1}\simeq E_l^{p-l,\,q+l-1}/\ker d_l^{p-l,\,q+l-1}$ for all $p,q,l$, if we set $p':=p-l$ and $q':=q+l-1$, we have $p'+q'=k-1$ when $p+q=k$ and the above isomorphism translates to

\begin{eqnarray}\nonumber E_r^k\simeq \bigoplus\limits_{p+q=k}E_\infty^{p,\,q} \oplus \bigoplus\limits_{l\geq r}\bigoplus\limits_{p'+q'=k-1}(E_l^{p',\,q'}/\ker d_l^{p',\,q'}) \oplus\bigoplus\limits_{l\geq r}\bigoplus\limits_{p+q=k}(E_l^{p,\,q}/\ker d_l^{p,\,q})\end{eqnarray}

\noindent for every $k$. Now, $\mbox{dim}\oplus_{p+q=k}E_\infty^{p,\,q} = b_k$ (the $k^{th}$ Betti number of $X$) thanks to (\ref{eqn:Froelicher_decomposition}), so taking dimensions in the above isomorphism, we get (\ref{eqn:dim_E_r^k}). 

\vspace{1ex}

$(ii)$\, Since $d_l^{p,\,q}:E_l^{p,\,q}/\ker d_l^{p,\,q}\longrightarrow \mbox{Im}\,d_l^{p,\,q}$ is an isomorphism of type $(l,\,-l+1)$ for all $l, p, q$, we get for all $l\geq r$:

$$d(L_r^{p,\,q})=\bigoplus\limits_{l\geq r}d_l^{p,\,q}(E_l^{p,\,q}/\ker d_l^{p,\,q}) \hspace{3ex} \mbox{and} \hspace{3ex} d_l^{p,\,q}(E_l^{p,\,q}/\ker d_l^{p,\,q})\subset E_l^{p+l,\,q-l+1}\subset C^\infty_{p+l,\,q-l+1}\subset{\cal A}_{p+r}^{p+q+1}$$

\noindent under the identification of each space $E_l^{p+l,\,q-l+1}$ with a subspace of $C^\infty_{p+l,\,q-l+1}$ defined by the isomorphisms (\ref{eqn:E_0_splitting}), (\ref{eqn:E_1_splitting}), (\ref{eqn:E_r_splitting}). This proves (\ref{eqn:d_L_r_pq}). \hfill $\Box$

\subsection{Explicit description of the above identifications}\label{subsection:explicit}

We take this opportunity to point out an explicit description of the differentials $d_r$ in cohomology and of their unique realisations induced by a given Hermitian metric on $X$.

\begin{Lem}\label{Lem:E_r-d_r_description} Let $X$ be a compact complex manifold with $\mbox{dim}_\C X=n$. 

\vspace{1ex}

$(i)$\, For every $r$ and every bidegree $(p,\,q)$, the vector space of type $(p,\,q)$ featuring on the $r^{th}$ page of the Fr\"olicher spectral sequence of $X$ can be explicitly described as the following set of multi-cohomology classes (i.e. each of these is the $d_{r-1}$-class of a $d_{r-2}$-class $\dots $ of a $d_1$-class of a $\bar\partial$-class): \begin{equation}\label{eqn:E_r_description}E_r^{p,\,q}=\{[\dots[[\alpha]_{\bar\partial}]_{d_1}\dots]_{d_{r-1}}\,\,\mid\,\,\alpha\in C^\infty_{p,\,q}(X,\,\C) \hspace{1ex} \mbox{such that} \hspace{1ex} \alpha \hspace{1ex} \mbox{satisfies condition} \hspace{1ex} (P_r)\},\end{equation} \noindent where condition $(P_r)$ on $\alpha$ requires the existence of forms $u_l\in C^\infty_{p+l,\,q-l}(X,\,\C)$ for $l\in\{1,\dots , r-1\}$ such that \begin{equation}\label{eqn:forms_u_l_existence}\bar\partial\alpha=0,\hspace{1ex}\partial\alpha=\bar\partial u_1,\hspace{1ex} \partial u_1 = \bar\partial u_2,\dots , \partial u_{r-2} = \bar\partial u_{r-1}.\end{equation}

$(ii)$\, For every $r$ and every bidegree $(p,\,q)$, the differential $d_r = d_r^{p,\,q}:E_r^{p,\,q}\longrightarrow E_r^{p+r,\,q-r+1}$ featuring on the $r^{th}$ page of the Fr\"olicher spectral sequence of $X$ is explicitly described as

\begin{equation}\label{eqn:d_r_description}d_r\bigg([\dots[[\alpha]_{\bar\partial}]_{d_1}\dots]_{d_{r-1}}\bigg) = [\dots[[\partial u_{r-1}]_{\bar\partial}]_{d_1}\dots]_{d_{r-1}}, \end{equation}

\noindent for every $[\dots[[\alpha]_{\bar\partial}]_{d_1}\dots]_{d_{r-1}}\in E_r^{p,\,q}$. Moreover, this description of $d_r$ is independent of the choice of forms $u_l\in C^\infty_{p+l,\,q-l}(X,\,\C)$ in (\ref{eqn:forms_u_l_existence}) (which are unique only modulo $\ker\bar\partial$).

\end{Lem}

\noindent {\it Proof.} These facts are well-known (cf. [CFGU97]). We will only explain the well-definedness of formula (\ref{eqn:d_r_description}) for $d_r$. Let $(u_1,\dots , u_{r-1})$ and $(u_1+\zeta_1,\dots , u_{r-1}+\zeta_{r-1})$ be two sets of forms satisfying (\ref{eqn:forms_u_l_existence}), i.e. $\bar\partial\alpha=0$, $\partial\alpha=\bar\partial u_1 = \bar\partial(u_1 + \zeta_1)$ and $$\partial u_1 = \bar\partial u_2 \hspace{1ex} \mbox{and} \hspace{1ex} \partial(u_1 + \zeta_1) = \bar\partial(u_2 + \zeta_2),\dots , \partial u_{r-2} = \bar\partial u_{r-1} \hspace{1ex} \mbox{and} \hspace{1ex} \partial(u_{r-2} + \zeta_{r-2}) = \bar\partial(u_{r-1} + \zeta_{r-1}).$$ \noindent These identities imply the identities $$\bar\partial\zeta_1 = 0, \hspace{1ex} \partial\zeta_1 = \bar\partial\zeta_2, \dots , \partial\zeta_{r-2} = \bar\partial\zeta_{r-1},$$ \noindent which, in turn, imply that $\zeta_1$ satisfies condition $(P_{r-1})$ (hence defines a multi-cohomology class lying in $E_{r-1}^{p+1,\,q-1}$) and that $$d_{r-1}([\dots[[\zeta_1]_{\bar\partial}]_{d_1}\dots]_{d_{r-2}}) = [\dots[[\partial\zeta_{r-1}]_{\bar\partial}]_{d_1}\dots]_{d_{r-2}}\in\mbox{Im}\,d_{r-1}.$$ \noindent Consequently, $[[\dots[[\partial\zeta_{r-1}]_{\bar\partial}]_{d_1}\dots]_{d_{r-2}}]_{d_{r-1}}=0$, so 

\vspace{1ex}

\hspace{15ex} $[\dots[[\partial(u_{r-1} + \zeta_{r-1})]_{\bar\partial}]_{d_1}\dots]_{d_{r-1}} = [\dots[[\partial u_{r-1}]_{\bar\partial}]_{d_1}\dots]_{d_{r-1}}.$ 

\vspace{1ex}

\noindent Thus, the result we get by formula (\ref{eqn:d_r_description}) for $d_r([\dots[[\alpha]_{\bar\partial}]_{d_1}\dots]_{d_{r-1}})$ is the same whether we work with the choices $(u_1,\dots , u_{r-1})$ or $(u_1+\zeta_1,\dots , u_{r-1}+\zeta_{r-1})$.   \hfill $\Box$

\vspace{3ex}

Thus, $d\alpha = \partial\alpha$ induces the multi-cohomology class $d_r([\dots[[\alpha]_{\bar\partial}]_{d_1}\dots]_{d_{r-1}})$. This helps to explain that, intuitively, $d$ acts as $d_r$ on representatives of $E_r$-classes (cf. (\ref{eqn:d_splitting})).

\vspace{3ex}

Now, recall that infinitely many choices are possible for the isomorphisms (\ref{eqn:E_0_splitting}), (\ref{eqn:E_1_splitting}) and (\ref{eqn:E_r_splitting}). However, any fixed Hermitian metric $\omega$ on $X$ selects a unique realisation of each of these isomorphisms and, implicitly, identifies each space $E_r^{p,\,q}$ with a precise subspace ${\cal H}_r^{p,\,q}$ (depending on $\omega$) of $C^\infty_{p,\,q}(X,\,\C)$ via an isomorphism $E_r^{p,\,q}\simeq{\cal H}_r^{p,\,q}$ depending on $\omega$. These {\it harmonic} subspaces ${\cal H}_r^{p,\,q}\subset C^\infty_{p,\,q}(X,\,\C)$ are constructed by induction on $r\geq 1$ as follows. 

\begin{Def}\label{Def:H-spaces_def} Let ${\cal H}_1^{p,\,q}\subset C^\infty_{p,\,q}(X,\,\C)$ be the orthogonal complement for the $L^2_\omega$-norm of $\mbox{Im}\,d_0^{p,\,q-1}$ in $\ker d_0^{p,\,q}$. Due to (\ref{eqn:E_0_splitting}), ${\cal H}_1^{p,\,q}$ is isomorphic to $E_1^{p,\,q}$. In every bidegree $(p,\,q)$, the linear map $d_1^{p,\,q}:E_1^{p,\,q}\longrightarrow E_1^{p+1,\,q}$ induces a linear map (denoted by the same symbol) $d_1^{p,\,q}:{\cal H}_1^{p,\,q}\longrightarrow {\cal H}_1^{p+1,\,q}$ via the isomorphisms ${\cal H}_1^{p,\,q}\simeq E_1^{p,\,q}$ and ${\cal H}_1^{p+1,\,q}\simeq E_1^{p+1,\,q}$. Let ${\cal H}_2^{p,\,q}\subset {\cal H}_1^{p,\,q}\subset C^\infty_{p,\,q}(X,\,\C)$ be the orthogonal complement for the $L^2_\omega$-norm of $\mbox{Im}\,d_1^{p-1,\,q}$ in $\ker d_1^{p,\,q}$ (viewed as subspaces of ${\cal H}_1^{p,\,q}$). Due to (\ref{eqn:E_1_splitting}), ${\cal H}_2^{p,\,q}$ is isomorphic to $E_2^{p,\,q}$. Continuing inductively, when the linear maps $d_r^{p,\,q}:E_r^{p,\,q}\longrightarrow E_r^{p+r,\,q-r+1}$ have induced counterparts (denoted by the same symbol) $d_r^{p,\,q}:{\cal H}_r^{p,\,q}\longrightarrow {\cal H}_r^{p+r,\,q-r+1}$ between the already constructed subspaces ${\cal H}_r^{p,\,q}\subset C^\infty_{p,\,q}(X,\,\C)$ and ${\cal H}_r^{p+r,\,q-r+1}\subset C^\infty_{p+r,\,q-r+1}(X,\,\C)$, we let ${\cal H}_{r+1}^{p,\,q}\subset {\cal H}_r^{p,\,q}\subset C^\infty_{p,\,q}(X,\,\C)$ be the orthogonal complement for the $L^2_\omega$-norm of $\mbox{Im}\,d_r^{p-r,\,q+r-1}$ in $\ker d_r^{p,\,q}$ (viewed as subspaces of ${\cal H}_r^{p,\,q}$). Due to (\ref{eqn:E_r_splitting}), ${\cal H}_{r+1}^{p,\,q}$ is isomorphic to $E_{r+1}^{p,\,q}$. 

\end{Def}

Note that we have \begin{eqnarray}\label{eqn:harmonic-spaces}\nonumber {\cal H}_1^{p,\,q} & = & \ker\,(\Delta'':C^\infty_{p,\,q}(X,\,\C)\longrightarrow C^\infty_{p,\,q}(X,\,\C)) = \{u\in C^\infty_{p,\,q}(X,\,\C)\,\mid\,\bar\partial u=0 \hspace{1ex} \mbox{and} \hspace{1ex} \bar\partial^\star u=0\}, \\
\nonumber {\cal H}_2^{p,\,q} & = & \ker\,(\widetilde\Delta:C^\infty_{p,\,q}(X,\,\C)\longrightarrow C^\infty_{p,\,q}(X,\,\C)) \\
 & = & \{u\in C^\infty_{p,\,q}(X,\,\C)\,\mid\,\bar\partial u=0, \hspace{1ex} \bar\partial^\star u=0, \hspace{1ex} p''(\partial u)=0 \hspace{1ex} \mbox{and} \hspace{1ex} p''\partial^\star u=0\},\end{eqnarray}

\noindent where $\widetilde\Delta=\partial p''\partial^\star + \partial^\star p''\partial + \Delta''$ is the pseudo-differential Laplacian constructed in [Pop16].

Also note that standard Hodge theory (for the elliptic differential operator $\Delta''$) is used to ensure that $\mbox{Im}\,d_0^{p,\,q-1}$ is closed in $C^\infty_{p,\,q}(X,\,\C)$ and that ${\cal H}_1^{p,\,q}$ is finite-dimensional. However, all the other images $\mbox{Im}\,d_r^{p-r,\,q+r-1}$ are automatically closed since they are (necessarily finite-dimensional) vector subspaces of a finite-dimensional vector space. It is also possible to construct pseudo-differential operators $\widetilde\Delta_{(r)}:C^\infty_{p,\,q}(X,\,\C)\longrightarrow C^\infty_{p,\,q}(X,\,\C)$ whose kernels are isomorphic to the spaces ${\cal H}_r^{p,\,q}$ (cf. forthcoming joint work of the author with L. Ugarte, where the Hodge theory found in [Pop16] for the second page of the Fr\"olicher spectral sequence is extended to all the pages), making these spaces into {\bf harmonic spaces} for these {\bf pseudo-differential Laplacians}, but the mere spaces ${\cal H}_r^{p,\,q}$ suffice for our purposes in this paper.

When the vector space $C^\infty_{p,\,q}(X,\,\C)$ is endowed with the $L^2$-norm induced by $\omega$, every subspace ${\cal H}_r^{p,\,q}$ inherits the restricted norm. On the other hand, every space $E_r^{p,\,q}$ has a quotient norm induced by the $L^2_\omega$-norm. The isomorphisms $E_r^{p,\,q}\simeq{\cal H}_r^{p,\,q}$ constructed above are isometries when $E_r^{p,\,q}$ and ${\cal H}_r^{p,\,q}$ are endowed with the quotient, resp. $L^2$ norms.

\begin{Conc}\label{Conc:E_r-d_r_metric-realisations} Let $X$ be a compact complex manifold and let $\omega$ be any Hermitian metric on $X$. Let $\dots\subset{\cal H}_{r+1}^{p,\,q}\subset{\cal H}_r^{p,\,q}\subset\dots\subset{\cal H}_1^{p,\,q}\subset C^\infty_{p,\,q}(X,\,\C)$ be the subspaces of Definition \ref{Def:H-spaces_def} induced by $\omega$.

 For every $r$ and every bidegree $(p,\,q)$, each class $[\dots[[\alpha]_{\bar\partial}]_{d_1}\dots]_{d_{r-1}}\in E_r^{p,\,q}$ contains a unique representative $\alpha\in{\cal H}_r^{p,\,q}$ (necessarily satisfying condition $(P_r)$). For $l\in\{1,\dots , r-1\}$, let $u_l\in C^\infty_{p+l,\,q-l}(X,\,\C)$ be the {\bf unique} solutions with {\bf minimal} $L^2_\omega$-norms of the equations

\begin{equation}\label{eqn:forms_u_l_existence_bis}\nonumber\bar\partial\alpha=0,\hspace{1ex}\partial\alpha=\bar\partial u_1,\hspace{1ex} \partial u_1 = \bar\partial u_2,\dots , \partial u_{r-2} = \bar\partial u_{r-1}\end{equation}

\noindent constructed inductively from one another. The well-known Neumann formula yields 

$$u_1=\Delta^{''-1}\bar\partial^\star(\partial\alpha) \hspace{2ex} \mbox{and} \hspace{2ex} u_l=\Delta^{''-1}\bar\partial^\star(\partial u_{l-1}) \hspace{2ex} \mbox{for} \hspace{2ex} l\in\{2,\dots , r-1\}.$$

\noindent In particular, the maps $\alpha\mapsto u_1$ and $u_{l-1}\mapsto u_l$ are linear.

For all $r,p,q$, we define the linear operator

$$T_r = T_r^{p,\,q}:{\cal H}_r^{p,\,q}\longrightarrow C^\infty_{p+r,\,q-r+1}(X,\,\C), \hspace{3ex} \alpha\mapsto T_r(\alpha):=\partial u_{r-1}.$$

\noindent Since ${\cal H}_r^{p,\,q}$ is finite-dimensional, $T_r$ is bounded, so there exists a constant $C_r^{p,\,q}>0$ such that

$$||T_r(\alpha)||=||\partial u_{r-1}||\leq C_r^{p,\,q}\,||\alpha|| \hspace{3ex} \mbox{for all} \hspace{1ex} \alpha\in{\cal H}_r^{p,\,q}.$$

\noindent It is easy to see that $T_r(\alpha)$ need not belong to ${\cal H}_r^{p+r,\,q-r+1}$ when $\alpha\in{\cal H}_r^{p,\,q}$. If we let $P_r^{p,\,q}: C^\infty_{p,\,q}(X,\,\C)\longrightarrow{\cal H}_r^{p,\,q}$ be the $L_\omega$-orthogonal projection onto ${\cal H}_r^{p,\,q}$, we get

$$||(P_r^{p,\,q}\circ T_r)(\alpha)||=||P_r^{p,\,q}(\partial u_{r-1})||\leq ||\partial u_{r-1}||\leq C_r^{p,\,q}\,||\alpha|| \hspace{3ex} \mbox{for all} \hspace{1ex} \alpha\in{\cal H}_r^{p,\,q}.$$

\end{Conc}

\section{Use of the rescaled Laplacians in the study of the Fr\"olicher spectral sequence}\label{section:Laplacians_Froelicher}

In this section, we prove Theorem \ref{The:main1}.

As in [ES89], [GS91], [ALK00], we consider the {\it spectrum distribution function} associated with any of the rescaled Laplacians $\Delta_h$, $\Delta_{\omega_h}$ in our context. Its definition and its study are made far simpler in this setting than in those references by the manifold $X$ being {\it compact} and by the Laplacians $\Delta'$, $\Delta''$ being {\it elliptic}. 

\begin{Def}\label{Def:spectrum-d-function} Let $(X,\,\omega)$ be a compact Hermitian manifold with $\mbox{dim}_\C X=n$. For every $k\in\{0,\dots , n\}$ and every constant $\lambda\geq 0$, let $N_h^k(\lambda)$ stand for the number of eigenvalues (counted with multiplicities) of $\Delta_h$ that are $\leq\lambda$.

\end{Def}

Replacing $\Delta_h$ with $\Delta_{\omega_h}$ does not change the spectrum distribution function $N_h^k:[0,\,+\infty)\longrightarrow\N$ since $\Delta_h$ and $\Delta_{\omega_h}$ have the same eigenvalues with the same multiplicities (cf. Corollary \ref{Cor:spectra}). Theorem \ref{The:main1} can be reworded as ensuring the existence of a constant $C>0$ independent of $h$ such that, for all $r$ and $k$, we have 

\begin{equation}\label{eqn:main_identity_rewording}\mbox{dim}\,E_r^k = N_h^k(Ch^{2r})  \hspace{2ex} \mbox{when} \hspace{1ex} 0<h\ll 1.\end{equation}

\subsection{The Efremov-Shubin variational principle}\label{subsection:variational}

The main technical ingredient we will need is the following variant of the {\it variational principle} proved in a more general context in [ES89] and used extensively thereafter (e.g. [GS91], [ALK00]) in settings different from ours. We adapt to our situation the result of [ES89]. 

\begin{Prop}(see e.g. Efremov-Shubin [ES89])\label{Prop:variational-principle} Let $(X,\,\omega)$ be a compact Hermitian manifold with $\mbox{dim}_\C X=n$. For every $k=0,\dots , 2n$ and every $\lambda\geq 0$, the following identity holds

\begin{equation}\label{eqn:variational-principle}N_h^k(\lambda) = F_h^{k-1}(\lambda) + b_k + F_h^k(\lambda),\end{equation}

\noindent where $b_k$ is the $k^{th}$ Betti number of $X$ and the function $F_h^k:[0,\,+\infty)\longrightarrow\N$ is defined by

\begin{equation}\label{eqn:F_hk_def}F_h^k(\lambda) = \sup\limits_L\mbox{dim}\,L,\end{equation}

\noindent where $L$ ranges over the {\bf closed} vector subspaces of the quotient space $C^\infty_k(X,\,\C)/\ker\,d$ on which the operator $d:C^\infty_k(X,\,\C)/\ker\,d\longrightarrow C^\infty_{k+1}(X,\,\C)$ induced by $d:C^\infty_k(X,\,\C)\longrightarrow C^\infty_{k+1}(X,\,\C)$ satisfies the following $L^2_{\omega_h}$-norm estimate:

\begin{equation}\label{eqn:d_omega_h_estimate}||d\zeta||_{\omega_h}\leq\sqrt{\lambda}\,||\zeta||_{\omega_h}, \hspace{3ex} \mbox{for every}\hspace{1ex} \zeta\in L.\end{equation}

\noindent (The understanding is that $||d\zeta||_{\omega_h}$ stands for the usual $L^2$-norm induced by the metric $\omega_h$, while $||\zeta||_{\omega_h}$ stands for the {\bf quotient norm} induced on $C^\infty_k(X,\,\C)/\ker\,d$ by the $L^2_{\omega_h}$-norm.)

\end{Prop}

We will present a detailed proof of this statement along the lines of [ES89] with a few minor simplifications afforded by our special setting where the manifold $X$ is {\bf compact} and the operator $\Delta_h$ is {\bf elliptic}. While a more general version for unbounded operators on $L^2$ spaces was needed in [ALK00], we stress that, in this context, we can confine ourselves to the case of operators on spaces of $C^\infty$ differential forms.

The main step is the following statement (a version of the classical Min-Max Principle) that was proved in a more general setting in [ES89].

\begin{Prop}\label{Prop:variational-principle_initial} Let $(X,\,\omega)$ be a {\bf compact} Hermitian manifold with $\mbox{dim}_\C X=n$. For an arbitrary $k\in\{0,\dots , 2n\}$, let $P:C^\infty_k(X,\,\C)\longrightarrow C^\infty_k(X,\,\C)$ be an {\bf elliptic}, self-adjoint and non-negative differential operator of order $\geq 1$. 

Then, for every $\lambda\geq 0$, the {\bf spectrum distribution function $N_k$ of $P$} (i.e. $N_k(\lambda)$ is defined to be the number of eigenvalues of $P$, counted with multiplicities, that are $\leq\lambda$) is given by the following identities (in which the suprema are actually maxima):

\begin{equation}\label{eqn:variational-principle_initial}N_k(\lambda) = \sup\limits_{L\in{\cal L}_\lambda^{(k)}}\mbox{dim}\,L = \sup\limits_{E\in{\cal P}_\lambda^{(k)}}\mbox{Tr}\,E,   \end{equation}

\noindent where ${\cal L}_\lambda^{(k)}$ stands for the set of {\bf closed} vector subspaces $L\subset C^\infty_k(X,\,\C)$ such that

$$\langle\langle Pu,\,u\rangle\rangle\leq\lambda||u||^2 \hspace{3ex} \mbox{for all}\hspace{1ex} u\in L,$$

\noindent while ${\cal P}_\lambda^{(k)}$ stands for the set of all bounded linear operators $E:C^\infty_k(X,\,\C)\longrightarrow C^\infty_k(X,\,\C)$ satisfying the conditions:

\vspace{1ex}

$(i)$\, $E^2=E=E^\star$ (i.e. $E$ is an {\bf orthogonal projection} w.r.t. the $L^2_\omega$ inner product);

\vspace{1ex}

$(ii)$\, $\langle\langle Pu,\,u\rangle\rangle\leq\lambda||u||^2$ for all $u\in\mbox{Im}\,E$.    

\vspace{1ex}

\noindent (In other words, $E$ is the orthogonal projection onto one of the subspaces $L\in{\cal L}_\lambda^{(k)}$, so $L=\mbox{Im}\,E$ for some $L\in{\cal L}_\lambda^{(k)}$.)

\end{Prop}

\noindent {\it Proof.} The second identity in (\ref{eqn:variational-principle_initial}) follows at once from the fact that the dimension of any closed subspace 
$L\subset C^\infty_k(X,\,\C)$ equals the trace of the orthogonal projection onto $L$. So, we only have to prove the first identity in (\ref{eqn:variational-principle_initial}).

Since $X$ is compact and $P$ is elliptic, self-adjoint and non-negative, the spectrum of $P$ is discrete and consists of non-negative eigenvalues, while there exists a countable orthonormal (w.r.t. the $L^2_\omega$-inner product) basis of $C^\infty_k(X,\,\C)$ (and of the Hilbert space $L^2_k(X,\,\C)$ of $L^2$ $k$-forms) consisting of eigenvectors of $P$. For every $\mu\geq 0$, let $E_P(\mu)\subset C^\infty_k(X,\,\C)$ be the eigenspace of $P$ corresponding to the eigenvalue $\mu$ (with the understanding that $E_P(\mu)=\{0\}$ if $\mu$ is not an actual eigenvalue). The spaces $E_P(\mu)$ are finite-dimensional and consist of $C^\infty$ forms since $P$ is assumed to be elliptic (hence also hypoelliptic) and $X$ is compact. 

 For every $\lambda\geq 0$, let $L_\lambda:=\bigoplus\limits_{0\leq\mu\leq\lambda}E_P(\mu)\subset C^\infty_k(X,\,\C)$. Thus, $L_\lambda$ is finite-dimensional and $\mbox{dim}\,L_\lambda = N_k(\lambda)$, while $\langle\langle Pu,\,u\rangle\rangle\leq\lambda\,||u||^2$ for all $u\in L_\lambda$. Hence $L_\lambda\in{\cal L}_\lambda^{(k)}$, so $N_k(\lambda) \leq \sup\limits_{L\in{\cal L}_\lambda^{(k)}}\mbox{dim}\,L$.

To prove the reverse inequality, let $\lambda\geq 0$ and let $L\in{\cal L}_\lambda^{(k)}$. The existence of an orthonormal basis of eigenvectors implies the orthogonal direct-sum decomposition

$$C^\infty_k(X,\,\C) = \bigoplus\limits_{0\leq\mu\leq\lambda}E_P(\mu)\oplus\bigoplus\limits_{\mu>\lambda}E_P(\mu).$$ 

\noindent In particular, $\oplus_{\mu>\lambda}E_P(\mu) = \ker E_\lambda$, where $E_\lambda$ is the orthogonal projection onto $\oplus_{0\leq\mu\leq\lambda}E_P(\mu)$.

Now, $\langle\langle Pu,\,u\rangle\rangle > \lambda||u||^2$ for all $u\in\oplus_{\mu>\lambda}E_P(\mu)\setminus\{0\}$, while $\langle\langle Pu,\,u\rangle\rangle\leq\lambda||u||^2$ for all $u\in L$. So, $L\cap\ker E_\lambda=L\cap\oplus_{\mu>\lambda}E_P(\mu) = \{0\}$. This implies that the restriction

$$E_{\lambda|L}:L\longrightarrow\mbox{Im}\,E_\lambda = \bigoplus\limits_{0\leq\mu\leq\lambda}E_P(\mu)$$

\noindent is {\bf injective}. In particular, $\mbox{dim}\,L\leq\mbox{dim}\,\oplus_{0\leq\mu\leq\lambda}E_P(\mu)=N_k(\lambda)$. Since $L$ has been chosen arbitrarily in ${\cal L}_\lambda^{(k)}$, we conclude that $\sup\limits_{L\in{\cal L}_\lambda^{(k)}}\mbox{dim}\,L\leq N_k(\lambda)$ and we are done.  \hfill $\Box$

\vspace{3ex}

The second step towards proving Proposition \ref{Prop:variational-principle} is the standard $3$-space decomposition used in Hodge theory. For every $k=0,\dots , 2n$, the operator $\Delta_{\omega_h}:C^\infty_k(X,\,\C)\longrightarrow C^\infty_k(X,\,\C)$ is elliptic and since the manifold $X$ is compact and $d^2=0$, we have the $L^2_{\omega_h}$-orthogonal decomposition: \begin{equation}\label{eqn:3-space_decomp}C^\infty_k(X,\,\C)={\cal H}^k_{\Delta_{\omega_h}}(X,\,\C)\oplus E_k(X,\,\C) \oplus E_k^\star(X,\,\C), \hspace{3ex} \mbox{where} \hspace{2ex} \ker d={\cal H}^k_{\Delta_{\omega_h}}(X,\,\C)\oplus E_k(X,\,\C),   \end{equation}

\noindent and where ${\cal H}^k_{\Delta_{\omega_h}}(X,\,\C)$ is the kernel of $\Delta_{\omega_h}:C^\infty_k(X,\,\C)\longrightarrow C^\infty_k(X,\,\C)$, $E_k(X,\,\C):=\mbox{Im}\,(d:C^\infty_{k-1}(X,\,\C)\longrightarrow C^\infty_k(X,\,\C))$ and $E_k^\star(X,\,\C):=\mbox{Im}\,(d^\star_{\omega_h}:C^\infty_{k+1}(X,\,\C)\longrightarrow C^\infty_k(X,\,\C))$ . 

Moreover, each of the three subspaces into which $C^\infty_k(X,\,\C)$ splits in (\ref{eqn:3-space_decomp}) is {\bf $\Delta_{\omega_h}$-invariant}, i.e.

$$\Delta_{\omega_h}({\cal H}^k_{\Delta_{\omega_h}}(X,\,\C))\subset{\cal H}^k_{\Delta_{\omega_h}}(X,\,\C), \hspace{2ex}  \Delta_{\omega_h}(E_k(X,\,\C))\subset E_k(X,\,\C), \hspace{2ex} \Delta_{\omega_h}(E_k^\star(X,\,\C))\subset E_k^\star(X,\,\C)$$ 
 
\noindent because $\Delta_{\omega_h}$ commutes with $d$ and with $d^\star_{\omega_h}$. The invariance implies that an $L^2_{\omega_h}$-orthonormal basis $\{e_i^k(h)\}_{i\in\N^\star}$ of $C^\infty_k(X,\,\C)$ consisting of eigenvectors for $\Delta_{\omega_h}$ (whose existence follows from the standard elliptic theory) can be chosen such that each $e_i^k(h)$ belongs to one and only one of the subspaces ${\cal H}^k_{\Delta_{\omega_h}}(X,\,\C)$, $E_k(X,\,\C)$ and $E_k^\star(X,\,\C)$. Let $0\leq\lambda_1^k(h)\leq\dots\leq\lambda_i^k(h)\leq\dots$ be the corresponding eigenvalues, counted with multiplicities, of the rescaled Laplacian $\Delta_h:C^\infty_k(X,\,\C)\longrightarrow C^\infty_k(X,\,\C)$ ($=$ those of $\Delta_{\omega_h}:C^\infty_k(X,\,\C)\longrightarrow C^\infty_k(X,\,\C)$). Thus, $\Delta_{\omega_h}e_i^k(h) =\lambda_i^k(h)\,e_i^k(h)$ for all $i$. 

Consequently, we can define functions $F_h^k:[0,\,+\infty)\longrightarrow\N$ and $G_h^k:[0,\,+\infty)\longrightarrow\N$ by $$F_h^k(\lambda):=\sharp\{i\,\mid\, e_i^k(h)\in E_k^\star(X,\,\C) \hspace{2ex}\mbox{and}\hspace{2ex} \lambda_i^k(h)\leq\lambda\}$$

\noindent and $$G_h^k(\lambda):=\sharp\{i\,\mid\, e_i^k(h)\in E_k(X,\,\C) \hspace{2ex}\mbox{and}\hspace{2ex} \lambda_i^k(h)\leq\lambda\}.$$

\noindent These definitions of $F_h^k$ and $G_h^k(\lambda)$ are independent of the choice of orthonormal basis $\{e_i^k(h)\}_{i\in\N^\star}$ of $C^\infty_k(X,\,\C)$ satisfying the above properties.

\begin{Lem}\label{Lem:description_F-G} The functions $F_h^k$ and $G_h^k$ are the spectrum distribution functions of the restrictions $\Delta_{\omega_h|E_k^\star(X,\,\C)}:E_k^\star(X,\,\C)\longrightarrow E_k^\star(X,\,\C)$, resp. $\Delta_{\omega_h|E_k(X,\,\C)}:E_k(X,\,\C)\longrightarrow E_k(X,\,\C)$. 

 In other words, they are described as follows: \begin{eqnarray}\label{eqn:description_F-G} F_h^k(\lambda) & = & \sup\limits_{L\in{\cal L}_\lambda^{''(k)}}\mbox{dim}\,L,  \\
\nonumber G_h^k(\lambda) & = & \sup\limits_{L\in{\cal L}_\lambda^{'(k)}}\mbox{dim}\,L      \end{eqnarray}

\noindent where ${\cal L}_\lambda^{''(k)}$ stands for the set of {\bf closed} vector subspaces $L\subset E^\star_k(X,\,\C)$ such that

\begin{equation}\label{eqn:estimate_F}||du||^2_{\omega_h}\leq\lambda||u||^2_{\omega_h} \hspace{3ex} \mbox{for all}\hspace{1ex} u\in L,\end{equation}

\noindent and ${\cal L}_\lambda^{'(k)}$ stands for the set of {\bf closed} vector subspaces $L\subset E_k(X,\,\C)$ such that

\begin{equation}\label{eqn:estimate_G}||d^\star_{\omega_h}u||^2_{\omega_h}\leq\lambda||u||^2_{\omega_h} \hspace{3ex} \mbox{for all}\hspace{1ex} u\in L.\end{equation}

\end{Lem}

 \noindent {\it Proof.} This is an immediate application of the {\it variational principle} of Proposition \ref{Prop:variational-principle_initial} to the restrictions $\Delta_{\omega_h|E_k^\star(X,\,\C)}:E_k^\star(X,\,\C))\longrightarrow E_k^\star(X,\,\C)$ and $\Delta_{\omega_h|E_k(X,\,\C)}:E_k(X,\,\C))\longrightarrow E_k(X,\,\C)$. Estimates (\ref{eqn:estimate_F}) and (\ref{eqn:estimate_G}) are consequences of the identity $\langle\langle\Delta_{\omega_h}u,\,u\rangle\rangle_{\omega_h} = ||du||^2_{\omega_h} + ||d^\star_{\omega_h}u||^2_{\omega_h}$ and of the fact that $d^\star_{\omega_h}u=0$ whenever $u\in E^\star_k(X,\,\C)$ (since $\mbox{Im}\,d^\star_{\omega_h}\subset\ker d^\star_{\omega_h}$) and that $du=0$ whenever $u\in E_k(X,\,\C)$ (since $\mbox{Im}\,d\subset\ker d$).   \hfill $\Box$

\vspace{2ex}

The last ingredient we need is the following very simple observation.

\begin{Lem}\label{Lem:F_k_G_k+1} For every $\lambda\geq 0$ and every $k\in\{-1, 0, \dots , 2n\}$, we have

$$F_h^k(\lambda) = G_h^{k+1}(\lambda) \hspace{3ex} \mbox{with the understanding that} \hspace{3ex}  F^{-1}_h(\lambda) = G^{2n+1}_h(\lambda)=0.$$

\end{Lem}

\noindent {\it Proof.} We know from the orthogonal decompositions (\ref{eqn:3-space_decomp}) that the restriction of $d$ to $E_k^\star(X,\,\C)$ is injective, so $$d_{|E_k^\star(X,\,\C)}:E_k^\star(X,\,\C)\longrightarrow E_{k+1}(X,\,\C)$$

\noindent is an isomorphism. Moreover, $d\Delta_{\omega_h} = \Delta_{\omega_h}d$, so whenever $\Delta_{\omega_h}u_i=\lambda_i^k(h)\,u_i$, we get $\Delta_{\omega_h}(du_i)=\lambda_i^k(h)\,(du_i)$. Combined with the above isomorphism, with the invariance of $E_k^\star(X,\,\C)$ under $\Delta_{\omega_h}$ and with the definitions of $F_k^h(\lambda)$ and $G_{k+1}^h(\lambda)$, this implies the contention.  \hfill $\Box$.

\vspace{3ex}

\noindent {\it Proof of Proposition \ref{Prop:variational-principle}.} Putting together (\ref{eqn:3-space_decomp}), the definitions of $F_h^k(\lambda)$ and $G_h^k(\lambda)$ and the fact that the Hodge isomorphism ${\cal H}^k_{\Delta_{\omega_h}}\simeq H_{DR}^k(X,\,\C)$ (which follows at once from (\ref{eqn:3-space_decomp})) implies $b_k=\mbox{dim}\,{\cal H}^k_{\Delta_{\omega_h}}$, we get

$$N_h^k(\lambda) = b_k + G_h^k(\lambda) + F_h^k(\lambda)$$ 

\noindent for all $k$ and all $\lambda\geq 0$. Using Lemma \ref{Lem:F_k_G_k+1}, this is equivalent to (\ref{eqn:variational-principle}).

On the other hand, the descriptions (\ref{eqn:description_F-G}) and (\ref{eqn:estimate_F}) of $F_h^k(\lambda)$ coincide with the descriptions (\ref{eqn:F_hk_def}) and (\ref{eqn:d_omega_h_estimate}) thanks to the isomorphism $E_k^\star(X,\,\C)\simeq C^\infty_k(X,\,\C)/\ker d$, which is another consequence of the decompositions (\ref{eqn:3-space_decomp}).  \hfill $\Box$

\subsection{Metric independence of asymptotics}\label{subsection:metric-independence}

Although the following statement has no impact on either the statement of Theorem \ref{The:main1} or its proof, we pause briefly to show, exactly as in the foliated case of [ALK00], that the asymptotics of the eigenvalues $\lambda_i^k(h)$ and of the spectrum distribution function $N_h^k$ as $h\downarrow 0$ depend only on the complex structure of $X$. The proof is an easy application of the Variational Principle of Proposition \ref{Prop:variational-principle}.

\begin{Prop}\label{Prop:CS_invariants} The asymptotics of the $\lambda_i^k(h)$'s and of $N_h^k$ as $h\downarrow 0$ are independent of the choice of Hermitian metric $\omega$.

\end{Prop}

\noindent {\it Proof.} We adapt to our setting the proof of the corresponding result in [ALK00]. Let $\omega$ and $\omega'$ be two Hermitian metrics on $X$. They induce, respectively, rescaled metrics $(\omega_h)_{h>0}$ and $(\omega'_h)_{h>0}$. Let $N_h^{'k}(\lambda) = F_h^{'k-1}(\lambda) + b_k + F_h^{'k}(\lambda)$ be the spectrum distribution function associated with the rescaled Laplacian $\Delta_{\omega'_h}:C^\infty_k(X,\,\C)\longrightarrow C^\infty_k(X,\,\C)$, written as in (\ref{eqn:variational-principle}). 

Since $X$ is compact, there exists a constant $C>0$ such that the respective $L^2$-norms satisfy the following inequalities in every bidegree $(p,\,q)$: $$\frac{1}{C}\,||\,\,||_\omega\leq||\,\,||_{\omega'}\leq C\,||\,\,||_\omega,  \hspace{3ex} \mbox{hence} \hspace{3ex} \frac{1}{C}\,||\,\,||_{\omega_h}\leq||\,\,||_{\omega'_h}\leq C\,||\,\,||_{\omega_h} \hspace{2ex} \mbox{on}\hspace{1ex} L^2_{p,\,q}(X,\,\C) \hspace{2ex} \mbox{for every}\hspace{1ex} h>0.$$

\noindent The constant $C$ is independent of $h>0$ thanks to Formula \ref{Formula:L2_comparison}. 

 Hence, for every $\zeta\in C^\infty_k(X,\,\C)/\ker\,d$ such that $||d\zeta||_{\omega_h}\leq\sqrt{\lambda}\,||\zeta||_{\omega_h}$, we get $||d\zeta||_{\omega'_h}\leq\sqrt{C^4\lambda}\,|\zeta||_{\omega'_h}$. Thanks to Proposition \ref{Prop:variational-principle}, this implies that

$$F_h^k(\lambda)\leq F_h^{'k}(C^4\lambda), \hspace{3ex} \lambda\geq 0, \hspace{1ex} h>0.$$

\noindent By symmetry, we also get $F_h^{'k}(\lambda)\leq F_h^k(C^4\lambda)$, so putting the last two inequalities together, we get

$$F_h^{'k}(C^{-4}\lambda)\leq F_h^k(\lambda)\leq F_h^{'k}(C^4\lambda), \hspace{3ex} \lambda\geq 0, \hspace{1ex} h>0.$$

\noindent The proof is complete.  \hfill $\Box$

\subsection{Proof of the inequality ``$\leq$'' in Theorem \ref{The:main1}}\label{subsection:proof_upper-bound_main}

We are now in a position to prove the following

\begin{The}\label{The:main_1st} Let $(X,\,\omega)$ be a compact Hermitian manifold with $\mbox{dim}_\C X=n$. For every $r$ and every $k=0,\dots , 2n$, the following inequality holds:

\begin{equation}\label{eqn:main_inequality1}\mbox{dim}\,E_r^k \leq \sharp\{i\,\mid\,\lambda_i^k(h)\in O(h^{2r}) \hspace{2ex} \mbox{as} \hspace{1ex} h\downarrow 0\}.\end{equation}

\end{The}

\noindent {\it Proof.} We have to prove the existence of a uniform constant $C>0$ such that $\mbox{dim}\,E_r^k \leq N_h^k(Ch^{2r})$ for all $r, k$ and all $0<h\ll 1$. Recall the following facts:

\vspace{1ex}

$(i)$\, $\mbox{dim}\,E_r^k=b_k + m_r^{k-1} + m_r^k$, where $m_r^k:=\mbox{dim}L_r^k$ and $L_r^k:=\bigoplus\limits_{p+q=k}L_r^{p,\,q} = \bigoplus\limits_{p+q=k}\bigoplus\limits_{l\geq r}(E_l^{p,\,q}/\ker d_l^{p,\,q})$

\vspace{1ex}

\hspace{20ex} (proved in (\ref{eqn:dim_E_r^k}) of Lemma \ref{Lem:E_r^k_dimension});

\vspace{1ex}

$(ii)$\, $N_h^k(\lambda) = b_k + F_h^{k-1}(\lambda) + F_h^k(\lambda)$ for all $\lambda\geq 0$

\vspace{1ex}

\hspace{20ex} (cf. (\ref{eqn:variational-principle}) of Proposition \ref{Prop:variational-principle}).

\noindent Thus, it suffices to prove that

\begin{equation}\label{eqn:sufficient_inequality1}m_r^k\leq F_h^k(Ch^{2r}) \hspace{3ex} \mbox{for all}\hspace{1ex} 0<h\ll 1,\end{equation}

\noindent for a uniform constant $C>0$ and for all $r$ and $k$. 

 Now, thanks to the definition (\ref{eqn:F_hk_def}) of $F_h^k$, to prove (\ref{eqn:sufficient_inequality1}) it suffices to prove that $L_r^k$ is one of the subspaces of $C^\infty_k(X,\,\C)/\ker d$ contributing to the definition of $F_h^k(Ch^{2r})$ for some uniform constant $C>0$. In other words, it suffices to prove that there exists $C>0$ such that

\begin{equation}\label{eqn:sufficient_inequality2}||d\zeta||_{\omega_h}\leq\sqrt{C}\,h^r\,||\zeta||_{\omega_h},  \hspace{3ex} \mbox{for all}\hspace{1ex} \zeta\in L_r^k \hspace{1ex} \mbox{and all}\hspace{1ex}    0<h\ll 1.\end{equation}

\noindent Meanwhile, every $\zeta\in L_r^k = \bigoplus\limits_{p+q=k}L_r^{p,\,q}$ splits uniquely as $\zeta=\sum_{p+q=k}\zeta^{p,\,q}$ with $\zeta^{p,\,q}\in L_r^{p,\,q}$ for all $p,q$. Thus, it suffices to prove that, for a uniform constant $C>0$, we have

\begin{equation}\label{eqn:sufficient_inequality2}||d\zeta^{p,\,q}||_{\omega_h}\leq\sqrt{C}\,h^r\,||\zeta^{p,\,q}||_{\omega_h},  \hspace{3ex} \mbox{for all}\hspace{1ex} p, q, \hspace{1ex} \mbox{all} \hspace{1ex} \zeta^{p,\,q}\in L_r^{p,\,q} \hspace{1ex} \mbox{and all}\hspace{1ex}    0<h\ll 1.\end{equation}

 This holds mainly because $d_r$ is of type $(r,\,-r+1)$, so $d_r$ increases the holomorphic degree by $r$ and thus the norm $|\,\,\,\,|_{\omega_h}$ brings out an extra factor $h^r$. Specifically, for every $\zeta^{p,\,q}\in L_r^{p,\,q}$, (\ref{eqn:d_L_r_pq}) of Lemma \ref{Lem:E_r^k_dimension} yields $d\zeta^{p,\,q}\in d(L_r^{p,\,q})\subset{\cal A}_{p+r}^{p+q-1}$. Therefore, the holomorphic degree of $d\zeta^{p,\,q}$ is $\geq p+r$, so from Formula \ref{Formula:L2_comparison} we get

$$||d\zeta^{p,\,q}||_{\omega_h}\leq\frac{h^{p+r}}{h^n}\,||d\zeta^{p,\,q}||_\omega \hspace{3ex} \mbox{for all} \hspace{1ex} p, q, \hspace{1ex} \mbox{all} \hspace{1ex} \zeta^{p,\,q}\in L_r^{p,\,q} \hspace{1ex} \mbox{and all} \hspace{1ex} 0<h<1.$$

\noindent Now, $L_r^{p,\,q}$ is a {\it finite-dimensional} vector subspace of $C^\infty_k(X,\,\C)/\ker d$, so there exists a constant $C_r>0$ (depending on $r, p, q$, but independent of $h$) such that $||d\zeta^{p,\,q}||_\omega\leq C_r\,||\zeta^{p,\,q}||_\omega$ for all $\zeta^{p,\,q}\in L_r^{p,\,q}$. Meanwhile, Formula \ref{Formula:L2_comparison} tells us again that $||\zeta^{p,\,q}||_\omega = (h^n/h^p)\,||\zeta^{p,\,q}||_{\omega_h}$, so putting the last three relations together, we get

$$||d\zeta^{p,\,q}||_{\omega_h}\leq C_r\,h^r\,||\zeta^{p,\,q}||_{\omega_h} \hspace{3ex} \mbox{for all} \hspace{1ex} p, q, \hspace{1ex} \mbox{all} \hspace{1ex} \zeta^{p,\,q}\in L_r^{p,\,q} \hspace{1ex} \mbox{and all} \hspace{1ex} 0<h<1.$$

\noindent This proves (\ref{eqn:sufficient_inequality2}) after setting $C:=\max_{0\leq r\leq N\atop 0\leq p,q\leq n}C_r^2>0$.

 The proof is complete. \hfill $\Box$

\vspace{3ex}

Note that $L_r^k$ is a vector space of classes of cohomology classes, rather than of differential forms, so what is meant by $L_r^k$ in the above proof is its image in $C^\infty_k(X,\,\C)/\ker\,d$ under the isometries explained in $\S.$\ref{subsection:explicit}. We can use these isometries, the identification of $d$ acting on ${\cal H}_r^{p,\,q}$ with $d_r$ and Conclusion \ref{Conc:E_r-d_r_metric-realisations} in the following way to make the above proof even more explicit. If we choose $\zeta^{p,\,q}$ to be the $\omega_h$-{\it harmonic} representative of its class (also denoted by $\zeta^{p,\,q}$) and to play the role of $\alpha$ of Conclusion \ref{Conc:E_r-d_r_metric-realisations}, we can re-write the above inequalities in a more detailed form as follows: \begin{eqnarray}\nonumber||d\zeta^{p,\,q}||_{\omega_h} & = & ||(P(\partial u_{r-1})||_{\omega_h} \leq \frac{h^{p+r}}{h^n}\,||(P\circ T)(\zeta^{p,\,q})||_\omega \\ 
\nonumber & \leq & \frac{h^{p+r}}{h^n}\,C_r\,||\zeta^{p,\,q}||_\omega = C_r\,h^r||\alpha||_{\omega_h},\end{eqnarray} \noindent where $P$ and $T$ are the linear maps $P_r^{p,\,q}$ and  $T_r^{p,\,q}$ (with indices removed) of Conclusion \ref{Conc:E_r-d_r_metric-realisations} that was used above, while $||\,\,||_{\omega_h}$ stands for the $L^2_{\omega_h}$-norm when applied to a form and for the induced quotient norm when applied to a class.

\subsection{Preliminaries to the proof of the inequality ``$\geq$'' in Theorem \ref{The:main1}}\label{subsection:preliminaries_proof_lower-bound_main}

We will need a few simple observations.

\begin{Lem}\label{Lem:inner-prod_Delta_comparison} Let $(X,\,\omega)$ be a compact Hermitian manifold with $\mbox{dim}_\C X=n$. For every bidegree $(p,\,q)$ and every $(p,\,q)$-form $u$ on $X$, the following identities hold:

\begin{equation}\label{eqn:inner-prod_Delta_comparison}\langle\langle\Delta_h u,\,u\rangle\rangle_\omega = h^{2(n-p)}\,\langle\langle\Delta_{\omega_h}u,\,u\rangle\rangle_{\omega_h} = h^{2(n-p)}\,(||du||^2_{\omega_h} + ||d^\star_{\omega_h}u||^2_{\omega_h}).\end{equation}

\end{Lem}

\noindent {\it Proof.} The latter identity is obvious, so we will only prove the former one. Since $u$ is of pure type, (\ref{eqn:pure-type_vanishing}) yields the first identity below, while the second identity follows from Formula \ref{Formula:L2_comparison}: \begin{eqnarray}\label{eqn:inner-prod_Delta_comparison_proof1}\nonumber\langle\langle\Delta_h u,\,u\rangle\rangle_\omega & = & h^2\,\langle\langle\Delta' u,\,u\rangle\rangle_\omega + \langle\langle\Delta'' u,\,u\rangle\rangle_\omega = h^2\,h^{2(n-p)}\,\langle\langle\Delta' u,\,u\rangle\rangle_{\omega_h} + h^{2(n-p)}\,\langle\langle\Delta'' u,\,u\rangle\rangle_{\omega_h}\\
\nonumber & = & h^{2(n-p)}\,\langle\langle\Delta_{\omega_h}u,\,u\rangle\rangle_{\omega_h}.\end{eqnarray}

\noindent The last identity followed again from (\ref{eqn:pure-type_vanishing}).  \hfill $\Box$

\begin{Lem}\label{Lem:d-splitting_estimates} Let $u\in C^\infty_{p,\,q}(X,\,\C)$ be an arbitrary form. Considering the splitting $d=d^{(k)} =\bigoplus\limits_{0\leq r\leq N-1\atop p+q=k}d_r^{p,\,q} :C^\infty_k(X,\,\C)\longrightarrow C^\infty_{k+1}(X,\,\C)$ of the operator $d$ (see (\ref{eqn:d_splitting})) and the splitting $$u=\sum\limits_{r=0}^{N-1}u_r + \ker\,d, \hspace{3ex} \mbox{implying} \hspace{3ex} du = \sum\limits_{r=0}^{N-1}d_ru_r,$$

\noindent with $u_r\in E_r^{p,\,q}/\ker d_r^{p,\,q}$ (see $\S.$\ref{section:differentials-Froelicher} and recall that $d_r:E_r^{p,\,q}/\ker d_r^{p,\,q}\longrightarrow \mbox{Im}\,d_r^{p,\,q}\subset C^\infty_{p+r,\,q-r+1}(X,\,\C)$ is an isomorphism), the following identity holds: \begin{equation}\label{eqn:d-splitting_estimates} h^{2(n-p)}\,||du||^2_{\omega_h} = \sum\limits_{r=0}^{N-1} h^{2r}\,||d_ru_r||^2_\omega \hspace{3ex} \mbox{for all} \hspace{1ex} h>0.\end{equation}

\end{Lem}

\noindent {\it Proof.} Since $d_r$ is of type $(r,\,-r+1)$, $d_ru_r$ is of type $(p+r,\,q-r+1)$, so the $d_ru_r$'s are mutually orthogonal (w.r.t. any metric) when $r$ varies. We get $$||du||^2_{\omega_h} = \sum\limits_{r=0}^{N-1} ||d_ru_r||^2_{\omega_h} = \sum\limits_{r=0}^{N-1} \frac{h^{2(p+r)}}{h^{2n}}\,||d_ru_r||^2_\omega,$$ 

\noindent where for the last identity we used Formula \ref{Formula:L2_comparison}.   \hfill $\Box$

\vspace{3ex}

\begin{Lem}\label{Lem:d_r_star_comparison} For every $r$ and every bidegree $(p,\,q)$, the formal adjoints of $d_r$ w.r.t. the metrics $\omega_h$ and $\omega$ compare as follows: \begin{equation}\label{eqn:d_r_star_comparison}(d_r)^\star_{\omega_h} = h^{2r}\,(d_r)^\star_\omega.\end{equation}

\noindent Consequently, for every form $u\in C^\infty_{p,\,q}(X,\,\C)$, the following counterpart of Lemma \ref{Lem:d-splitting_estimates} for the adjoints holds. Considering the splitting $(d^{(k)})^\star_{\omega_h} =\bigoplus\limits_{0\leq r\leq N-1\atop p+q=k}(d_r^{p,\,q})^\star_{\omega_h} :C^\infty_{k+1}(X,\,\C)\longrightarrow C^\infty_k(X,\,\C)$ of the operator $d^\star$ and the splitting $$u=\sum\limits_{r=0}^{N-1}v_r + \ker\,d^\star_{\omega_h}, \hspace{3ex} \mbox{implying} \hspace{3ex} d^\star_{\omega_h} u = \sum\limits_{r=0}^{N-1}(d_r)^\star_{\omega_h}v_r,$$

\noindent with $v_r\in\mbox{Im}\, d_r^{p-r,\,q+r-1}$ (see $\S.$\ref{subsection:d_r-d_identification}), the following identity holds: \begin{equation}\label{eqn:d-splitting_estimates_adjoints} h^{2(n-p)}\,||d^\star_{\omega_h}u||^2_{\omega_h} = \sum\limits_{r=0}^{N-1} h^{2r}\,||(d_r)^\star_\omega v_r||^2_\omega \hspace{3ex} \mbox{for all} \hspace{1ex} h>0.\end{equation}

\end{Lem}

\noindent {\it Proof.} For every $(p,\,q)$-form $v$ and every $(p-r,\,q+r-1)$-form $u$, we have \begin{eqnarray}\nonumber\frac{h^{2(p-r)}}{h^{2n}}\, \langle\langle(d_r)^\star_{\omega_h}v,\,u\rangle\rangle_\omega    =\langle\langle(d_r)^\star_{\omega_h}v,\,u\rangle\rangle_{\omega_h} = \langle\langle v,\,d_ru\rangle\rangle_{\omega_h} = \frac{h^{2p}}{h^{2n}}\,\langle\langle v,\,d_ru\rangle\rangle_\omega = \frac{h^{2p}}{h^{2n}}\,\langle\langle (d_r)^\star_\omega v,\,u\rangle\rangle_\omega.\end{eqnarray}

\noindent This proves (\ref{eqn:d_r_star_comparison}). Using the mutual orthogonality of the $(d_r)^\star_{\omega_h}v_r$'s (due to bidegree reasons) and Formula \ref{Formula:L2_comparison}, we get $$||d^\star_{\omega_h}u||^2_{\omega_h} = \sum\limits_{r=0}^{N-1}||(d_r)^\star_{\omega_h}v_r||^2_{\omega_h} = \sum\limits_{r=0}^{N-1}\frac{h^{2(p-r)}}{h^{2n}}\,||(d_r)^\star_{\omega_h}v_r||^2_\omega = \sum\limits_{r=0}^{N-1}\frac{h^{2(p-r)}}{h^{2n}}\,h^{4r}\,||(d_r)^\star_\omega v_r||^2_\omega.$$

\noindent This proves (\ref{eqn:d-splitting_estimates_adjoints}).   \hfill $\Box$

\vspace{3ex}

 Putting together (\ref{eqn:inner-prod_Delta_comparison}), (\ref{eqn:d-splitting_estimates}) and (\ref{eqn:d-splitting_estimates_adjoints}), we get

\begin{Cor}\label{Cor:Delta_huu_sum_degrees-h} Let $(X,\,\omega)$ be a compact Hermitian manifold with $\mbox{dim}_\C X=n$. For every bidegree $(p,\,q)$ and every $(p,\,q)$-form $u$ on $X$, the following identity holds:

$$\langle\langle\Delta_h u,\,u\rangle\rangle_\omega = \sum\limits_{r'=0}^{N-1}h^{2r'}\,||d_{r'}u_{r'}||_\omega^2 + \sum\limits_{r'=0}^{N-1}h^{2r'}\,||(d_{r'})^\star_\omega v_{r'}||_\omega^2,$$

\noindent where $u$ splits uniquely (cf. $\S.$\ref{subsection:d_r-d_identification}) as $$u=\sum\limits_{r'=0}^{N-1}u_{r'} + \ker\,d = \sum\limits_{r'=0}^{N-1}v_{r'} + \ker\,d^\star = \sum\limits_{r'=0}^{N-1}u_{r'} + \sum\limits_{r'=0}^{N-1}v_{r'} + w$$

\noindent with $u_{r'}\in E_{r'}^{p,\,q}/\ker d_{r'}^{p,\,q}$, $v_{r'}\in\mbox{Im}\, d_{r'}^{p-r',\,q+r'-1}$ and $w\in E_\infty^{p,\,q}$.

\end{Cor}

\subsection{Proof of the inequality ``$\geq$'' in Theorem \ref{The:main1}}\label{subsection:proof_lower-bound_main}

Following again the analogy with the foliated case of [ALK00], we will actually prove a stronger statement from which the following result will follow as a corollary.

\begin{The}\label{The:main_2nd} Let $(X,\,\omega)$ be a compact Hermitian manifold with $\mbox{dim}_\C X=n$. For every $r$ and every $k=0,\dots , 2n$, the following inequality holds:

\begin{equation}\label{eqn:main_inequality2}\mbox{dim}\,E_r^k \geq \sharp\{i\,\mid\,\lambda_i^k(h)\in O(h^{2r}) \hspace{2ex} \mbox{as} \hspace{1ex} h\downarrow 0\}.\end{equation}

\end{The}

The first main ingredient we will use is the {\it pseudo-differential Laplacian}

$$\widetilde\Delta = \partial p''\partial^\star + \partial^\star p''\partial + \Delta'':C^\infty_{p,\,q}(X,\,\C)\longrightarrow C^\infty_{p,\,q}(X,\,\C)$$

\noindent defined in arbitrary bidegree $(p,\,q)$ and introduced in [Pop16], where $p'': C^\infty_{p,\,q}(X,\,\C)\longrightarrow\ker\Delta''$ is the orthogonal projection (w.r.t. the $L^2_\omega$-norm) onto the $\Delta''$-harmonic subspace of $C^\infty_{p,\,q}(X,\,\C)$. The pseudo-differential Laplacian $\widetilde\Delta$ gives a Hodge theory for the second page of the Fr\"olicher spectral sequence in the sense that there is a {\it Hodge isomorphism} \begin{equation}\label{eqn:Hodge_isom_pseudo-diff}E_2^{p,\,q}\stackrel{\simeq}{\longrightarrow}{\cal H}^{p,\,q}_{\widetilde\Delta}(X,\,\C):=\ker(\widetilde\Delta:C^\infty_{p,\,q}(X,\,\C)\longrightarrow C^\infty_{p,\,q}(X,\,\C)) \hspace{3ex} \mbox{for all} \hspace{1ex} p,q=0, \dots , n.\end{equation}

\noindent Note that $(p'')^2=p''=(p'')^\star$, so $\partial p''\partial^\star = (p''\partial^\star)^\star(p''\partial^\star)$ and $\partial^\star p''\partial=(p''\partial)^\star(p''\partial)$. Thus, $\widetilde\Delta$ is a sum of non-negative operators, so its kernel is the intersection of the respective kernels. Since $\ker(A^\star A)= \ker A$ for any operator $A$, we get

$$\ker\widetilde\Delta = \ker(p''\partial)\cap\ker(p''\partial^\star)\cap\ker\bar\partial\cap\ker\bar\partial^\star.$$

 The second main ingredient we will use is the following lower estimate of the rescaled Laplacian $\Delta_h$. It is the analogue in our context of a result in [ALK00].

\begin{Lem}\label{Lem:Delta_h_lower-estimate} Let $(X,\,\omega)$ be a compact Hermitian manifold with $\mbox{dim}_\C X=n$. There exists a constant $C>0$ such that the following inequality of linear operators (cf. Notation \ref{Not:introd}) holds on differential forms of any degree $k=0,\dots , 2n$:

$$\Delta_h\geq \frac{3}{4}\,\Delta'' + h^2\,\Delta' - Ch^2 \hspace{3ex} \mbox{for all} \hspace{1ex} h>0,$$

\noindent where $\Delta''=\bar\partial\bar\partial^\star + \bar\partial^\star\bar\partial$ and $\Delta'=\partial\partial^\star + \partial^\star\partial$ are the usual $\bar\partial$- and $\partial$-Laplacians.

\end{Lem} 

 The coefficients $3/4$ and $1$ are not optimal, but they suffice for our purposes and the proof provided below shows that they can be made optimal if this is desired.

\vspace{3ex}

\noindent {\it Proof of Lemma \ref{Lem:Delta_h_lower-estimate}.} We know from $(ii)$ of Lemma \ref{Lem:relations_Laplacians} that \begin{eqnarray}\label{eqn:Delta_h_bracket-expression}\Delta_h = \Delta'' + h^2\Delta' - h([\tau,\,\bar\partial^\star] + [\tau^\star,\,\bar\partial]),\end{eqnarray}

\noindent where $\tau=\tau_\omega:=[\Lambda,\,\partial\omega\wedge\cdot]$ is the zero-th order {\it torsion operator} of type $(1,\,0)$ associated with $\omega$.

 For any form $u$, the first-order terms on the r.h.s. of (\ref{eqn:Delta_h_bracket-expression}) are easily estimated using the Cauchy-Schwarz inequality as follows: \begin{eqnarray}\nonumber h\,|\langle\langle[\tau,\,\bar\partial^\star]u + [\tau^\star,\,\bar\partial]u,\,u\rangle\rangle| & = &  h\,|\langle\langle\bar\partial^\star u,\,\tau^\star u\rangle\rangle + \langle\langle\tau u,\,\bar\partial u\rangle\rangle + \langle\langle\bar\partial u,\,\tau u\rangle\rangle + \langle\langle\tau^\star u,\,\bar\partial^\star u\rangle\rangle| \\
\nonumber & \leq & 2h||\tau u||\,||\bar\partial u|| + 2h||\tau^\star u||\,||\bar\partial^\star u|| \\
\nonumber & \leq & \frac{1}{4}\,(||\bar\partial u||^2 + ||\bar\partial^\star u||^2) + 4h^2\,(||\tau u||^2 + ||\tau^\star u||^2) \\
\nonumber &\leq & \frac{1}{4}\, \langle\langle\Delta'' u,\, u\rangle\rangle + C h^2\,||u||^2,\end{eqnarray}

\noindent where the constant $C>0$ exists because the linear operators $\tau$ and $\tau^\star$ are of order zero, hence bounded. In particular, we get the operator inequality $- h([\tau,\,\bar\partial^\star] + [\tau^\star,\,\bar\partial])\geq -\frac{1}{4}\,\Delta'' - C h^2$ which, alongside (\ref{eqn:Delta_h_bracket-expression}), proves the contention.  \hfill $\Box$

\vspace{3ex}  

We are now ready to state and prove a general result that will imply Theorem \ref{The:main_2nd}.

\begin{The}\label{The:eigenforms_convergence} Let $(X,\,\omega)$ be a compact Hermitian manifold with $\mbox{dim}_\C X=n$. Let $k\in\{0,\dots, 2n\}$ and $r\geq 1$ be fixed integers. Suppose there exist a sequence $(h_i)_{i\in\N}$ of constants $h_i>0$ such that $h_i\downarrow 0$ and a sequence $(u_i)_{i\in\N}$ of $k$-forms $u_i\in C^\infty_k(X,\,\C)$ such that $||u_i||_\omega =1$ for every $i$ and 

\begin{equation}\label{eqn:hypothesis_H}\langle\langle\Delta_{h_i}u_i,\,u_i\rangle\rangle_\omega\in o(h_i^{2(r-1)}) \hspace{3ex} \mbox{as} \hspace{1ex} i\rightarrow +\infty.\end{equation}

\noindent Then, there exists a subsequence $(u_{i_l})_{l\in\N}$ of $(u_i)_{i\in\N}$ such that $(u_{i_l})_{l\in\N}$ converges in the $L^2_\omega$-topology to some $k$-form $u\in{\cal H}_r^k:=\oplus_{p+q=k}{\cal H}_r^{p,\,q}\simeq E_r^k$, where the ${\cal H}_r^{p,\,q}\subset C^\infty_{p,\,q}(X,\,\C)$ are the ``harmonic'' vector subspaces of Definition \ref{Def:H-spaces_def} induced by the metric $\omega$.

\end{The}

\noindent {\it Proof.} $\bullet$ {\bf Case $r=1$.} In this case, Hypothesis (\ref{eqn:hypothesis_H}) means that $\langle\langle\Delta_{h_i}u_i,\,u_i\rangle\rangle_\omega\longrightarrow 0$ as $i\rightarrow +\infty$. Then also $\langle\langle\Delta_{h_i}u_i,\,u_i\rangle\rangle_\omega + Ch_i^2\longrightarrow 0$ as $i\rightarrow +\infty$. Since, by Lemma \ref{Lem:Delta_h_lower-estimate}, we have

$$\langle\langle\Delta_{h_i}u_i,\,u_i\rangle\rangle_\omega + Ch_i^2\geq\frac{3}{4}\,\langle\langle\Delta''u_i,\,u_i\rangle\rangle_\omega + h_i^2\,\langle\langle\Delta'u_i,\,u_i\rangle\rangle_\omega\geq 0 \hspace{3ex} \mbox{for all}\hspace{1ex} i\in\N,$$

\noindent we get \begin{equation}\label{eqn:1st_convergence}(i)\,\,\,\langle\langle\Delta''u_i,\,u_i\rangle\rangle_\omega\longrightarrow 0 \hspace{3ex} \mbox{as}\hspace{1ex} i\rightarrow +\infty \hspace{3ex} \mbox{and} \hspace{3ex} (ii)\,\,\,h_i^2\,\langle\langle\Delta'u_i,\,u_i\rangle\rangle_\omega\longrightarrow 0 \hspace{3ex} \mbox{as}\hspace{1ex} i\rightarrow +\infty.\end{equation}

 Meanwhile, the $\bar\partial$-Laplacian $\Delta''$ is {\it elliptic} and the manifold $X$ is {\it compact}, so the G$\mathring{a}$rding inequality yields constants $\delta_1,\delta_2>0$ such that the first inequality below holds:

$$\delta_2\,||u_i||_{W^1}\leq\langle\langle\Delta''u_i,\,u_i\rangle\rangle_\omega + \delta_1\,||u_i||_\omega\leq C_1, \hspace{3ex} \mbox{for all}\hspace{1ex} i\in\N,$$

\noindent where $||\,\,\,||_{W^1}$ stands for the Sobolev norm $W^1$ induced by the metric $\omega$. The second inequality above holds for some constant $C_1>0$ since the quantity $\langle\langle\Delta''u_i,\,u_i\rangle\rangle_\omega$ converges to zero (cf. (\ref{eqn:1st_convergence})), hence is bounded, and $||u_i||_\omega = 1$ by the hypothesis of Theorem \ref{The:eigenforms_convergence}.

 Consequently, the sequence $(u_i)_{i\in\N}$ is bounded in the Sobolev space $W^1$ (a Hilbert space), so by the Banach-Alaoglu Theorem there exists a subsequence $(u_{i_l})_{l\in\N}$ that converges in the weak topology of $W^1$ to some $k$-form $u\in W^1$. In particular, the following convergences hold in the weak topology of distributions: $$\bar\partial u_{i_l}\longrightarrow \bar\partial u  \hspace{3ex} \mbox{and} \hspace{3ex} \bar\partial^\star u_{i_l}\longrightarrow \bar\partial^\star u \hspace{3ex} \mbox{as} \hspace{3ex} l\rightarrow +\infty.$$

\noindent On the other hand, $||\bar\partial u_i||^2 + ||\bar\partial^\star u_i||^2 = \langle\langle\Delta''u_i,\,u_i\rangle\rangle_\omega\longrightarrow 0$ as $i\rightarrow +\infty$, so $\bar\partial u_i\longrightarrow 0$ and $\bar\partial^\star u_i\longrightarrow 0$ in the $L^2$-topology as $i\rightarrow +\infty$. Comparing this with the above convergences in the weak topology of distributions, we get $$\bar\partial u = 0 \hspace{3ex} \mbox{and} \hspace{3ex} \bar\partial^\star u=0,$$

\noindent which, by (\ref{eqn:harmonic-spaces}), is equivalent to $u\in\ker\,(\Delta'':C^\infty_k(X,\,\C)\longrightarrow C^\infty_k(X,\,\C))={\cal H}_1^k\simeq E_1^k$.

Note that by the Rellich Lemma (asserting the compactness of the inclusion $W^1\hookrightarrow L^2$), the convergence of $(u_{i_l})_{l\in\N}$ to $u$ in the weak topology of $W^1$ implies that $(u_{i_l})_{l\in\N}$ also converges in the $L^2$-topology to $u$. Moreover, the ellipticity of $\Delta''$ and the relation $u\in\ker\Delta''$ imply that $u$ is $C^\infty$.  

\vspace{1ex}

$\bullet$ {\bf Case $r=2$.} In this case, Hypothesis (\ref{eqn:hypothesis_H}) means that $\langle\langle\Delta_{h_i}u_i,\,u_i\rangle\rangle_\omega\in o(h_i^2)$ as $i\rightarrow +\infty$. Since $\langle\langle\Delta_{h_i}u_i,\,u_i\rangle\rangle_\omega = ||d_{h_i}u_i||^2 + ||d_{h_i}^\star u_i||^2 = ||h_i\partial u_i + \bar\partial u_i||^2 + ||h_i\partial^\star u_i + \bar\partial^\star u_i||^2$, this implies that \begin{equation}\label{eqn:2nd_convergence}\partial u_i + \frac{1}{h_i}\,\bar\partial u_i\longrightarrow 0 \hspace{2ex} \mbox{and} \hspace{2ex} \partial^\star u_i + \frac{1}{h_i}\,\bar\partial^\star u_i\longrightarrow 0  \hspace{3ex} \mbox{in the} \hspace{1ex} L^2\mbox{-topology, as} \hspace{1ex} i\rightarrow +\infty.\end{equation}

\noindent Since the orthogonal projection $p''$ onto $\ker\Delta''$ is continuous w.r.t. the $L^2$-topology and since $p''\bar\partial =0$ and $p''\bar\partial^\star=0$ (because $\mbox{Im}\,\bar\partial\perp\ker\Delta''$ and $\mbox{Im}\,\bar\partial^\star\perp\ker\Delta''$), an application of $p''$ to (\ref{eqn:2nd_convergence}) yields

\begin{equation}\label{eqn:3nd_convergence}p''\partial u_i\longrightarrow 0 \hspace{2ex} \mbox{and} \hspace{2ex} p''\partial^\star u_i\longrightarrow 0  \hspace{3ex} \mbox{in the} \hspace{1ex} L^2\mbox{-topology, as} \hspace{1ex} i\rightarrow +\infty.\end{equation}

 On the other hand, we know from the discussion of the case $r=1$ (whose weaker assumption is still valid in the case $r=2$) that there exists a subsequence $(u_{i_l})_{l\in\N}$ that converges in the weak topology of $W^1$ to some $k$-form $u\in W^1$. Thus, $\partial u_{i_l}\longrightarrow\partial u\in L^2$ in the weak topology of $L^2$ as $l\rightarrow +\infty$. This means that $$\langle\langle\partial u_{i_l},\,v\rangle\rangle_\omega\longrightarrow\langle\langle\partial u,\,v\rangle\rangle_\omega \hspace{3ex} \mbox{for all} \hspace{1ex} v\in L^2,  \hspace{3ex} \mbox{hence} \hspace{1ex} \langle\langle\partial u_{i_l},\,p''v\rangle\rangle_\omega\longrightarrow\langle\langle\partial u,\,p''v\rangle\rangle_\omega \hspace{3ex} \mbox{for all} \hspace{1ex} v\in L^2,$$

\noindent as $l\rightarrow +\infty$. (The second convergence follows from the first since $||p''v||\leq ||v||$ for all $v\in L^2$, so $p''(L^2)\subset L^2$.) Now, $p''$ is self-adjoint, so the last convergence translates to

$$\langle\langle p''\partial u_{i_l},\,v\rangle\rangle_\omega\longrightarrow\langle\langle p''\partial u,\,v\rangle\rangle_\omega \hspace{3ex} \mbox{as}\hspace{1ex} l\rightarrow +\infty, \hspace{1ex} \mbox{for all} \hspace{1ex} v\in L^2.$$

 This means that $p''\partial u_{i_l}$ converges to $p''\partial u$ in the weak topology of $L^2$ as $l\rightarrow +\infty$. However, we know from (\ref{eqn:3nd_convergence}) that $p''\partial u_{i_l}$ converges to $0$ in the $L^2$-topology. Hence $p''\partial u =0$. The same argument run with $\partial^\star$ in place of $\partial$ yields that $p''\partial^\star u =0$. On the other hand, we know from the discussion of the case $r=1$ that $u\in\ker\bar\partial\cap\ker\bar\partial^\star=\ker\Delta''$, so we get

$$u\in\ker(p''\partial)\cap\ker(p''\partial^\star)\cap\ker\bar\partial\cap\ker\bar\partial^\star = {\cal H}_2^k \simeq E_2^k$$

\noindent after remembering the description (\ref{eqn:harmonic-spaces}) of the spaces ${\cal H}_2^{p,\,q}$ and that ${\cal H}_2^k=\oplus_{p+q=k}{\cal H}_2^{p,\,q} $.

\vspace{1ex}

$\bullet$ {\bf Case $r\geq 3$.} Using the information from the first two cases and from subsection $\S.$\ref{subsection:preliminaries_proof_lower-bound_main}, this last case can easily be dealt with as follows. 

For each of the $k$-forms $u_i$ given by the hypotheses of Theorem \ref{The:eigenforms_convergence}, we consider the splitting $$u_i= \sum\limits_{r'=0}^{N-1}u^{(i)}_{r'} + \sum\limits_{r'=0}^{N-1}v^{(i)}_{r'} + w_i,$$  \noindent with $u^{(i)}_{r'}\in E_{r'}^{p,\,q}/\ker d_{r'}^{p,\,q}$, $v^{(i)}_{r'}\in\mbox{Im}\, d_{r'}^{p-r',\,q+r'-1}$ and $w_i\in E_\infty^{p,\,q}$, and the corresponding splitting $$\langle\langle\Delta_{h_i} u_i,\,u_i\rangle\rangle_\omega = \sum\limits_{r'=0}^{N-1}h_i^{2r'}\,||d_{r'}u^{(i)}_{r'}||_\omega^2 + \sum\limits_{r'=0}^{N-1}h_i^{2r'}\,||(d_{r'})^\star_\omega v^{(i)}_{r'}||_\omega^2$$

\noindent obtained in Corollary \ref{Cor:Delta_huu_sum_degrees-h}

 On the other hand, (\ref{eqn:hypothesis_H}) ensures that $\langle\langle\Delta_{h_i} u_i,\,u_i\rangle\rangle_\omega\in o(h_i^{2(r-1)})$ as $i\rightarrow +\infty$. Together with the above identity, this implies the following convergences in the $L^2_\omega$-norm as $i\rightarrow +\infty$: $$d_{r'}u^{(i)}_{r'}\longrightarrow 0  \hspace{3ex} \mbox{and} \hspace{3ex} (d_{r'})^\star_\omega v^{(i)}_{r'}\longrightarrow 0 \hspace{3ex} \mbox{for every} \hspace{1.5ex} r'\in\{0,\dots , r-1\}.$$ 

\noindent We even get $$\frac{1}{h_i^{r-r'-1}}\,d_{r'}u^{(i)}_{r'}\longrightarrow 0  \hspace{3ex} \mbox{and} \hspace{3ex} \frac{1}{h_i^{r-r'-1}}\,(d_{r'})^\star_\omega v^{(i)}_{r'}\longrightarrow 0 \hspace{3ex} \mbox{for every} \hspace{1.5ex} r'\in\{0,\dots , r-1\}.$$

\noindent Defining in an ad hoc way a ``formal'' Laplacian by $\Delta_{r'}^{formal}:=d_{r'}(d_{r'})_\omega^\star + (d_{r'})_\omega^\star d_{r'}$, we get that the limit $u$ of a subsequence of $(u_i)_{i\in\N}$ lies in $$\ker\bigg(\Delta_{r-1}^{formal}:\bigoplus\limits_{p+q=k}E_{r-1}^{p,\,q}\longrightarrow\bigoplus\limits_{p+q=k}E_{r-1}^{p,\,q}\bigg)\simeq{\cal H}_r^k\simeq E_r^k$$

\noindent and we are done. \hfill $\Box$

\vspace{3ex}

\noindent {\it Proof of Theorem \ref{The:main_2nd}.} It is an immediate consequence of Theorem \ref{The:eigenforms_convergence}. Indeed, fix any $r\in\N^\star$ and $k\in\{0,\dots , 2n\} $ and suppose that inequality (\ref{eqn:main_inequality2}) does not hold. Then, the reverse strict inequality holds, so there exists a sequence $(h_i)_{i\in\N}$ of positive constants such that $h_i\downarrow 0$ when $i\rightarrow +\infty$ and a sequence $(u_i)_{i\in\N}$ of eigenvectors for the Laplacians $\Delta_{h_i}$ acting on $k$-forms such that $||u_i||_\omega=1$, $u_i\perp{\cal H}_r^k$ for all $i$ and $\langle\langle\Delta_{h_i}u_i,\,u_i\rangle\rangle\in o(h_i^{2(r-1)})$ as $i\rightarrow +\infty$.

Thanks to Theorem \ref{The:eigenforms_convergence}, there exists a subsequence $(u_{i_l})_{l\in\N}$ of $(u_i)_{i\in\N}$ such that $(u_{i_l})_{l\in\N}$ converges in the $L^2_\omega$-topology to some $k$-form $u\in{\cal H}_r^k\simeq E_r^k$. However, the form $u$ is orthogonal to ${\cal H}_r^k$ since $u_i\perp{\cal H}_r^k$ for all $i$ and the orthogonality property is preserved in the limit. Since $||u||_\omega=1$ (because $||u_i||_\omega=1$ for all $i$), $u\neq 0$, so $u$ cannot be at once orthogonal to and a member of ${\cal H}_r^k$. This is a contradiction. \hfill $\Box$

\section{Consequences of Theorem \ref{The:main1}}\label{section:consequences_MT1}

The following consequences of Theorem \ref{The:main1} are of independent interest.

\begin{Prop}\label{Prop:numerical-F-duality} Let $X$ be a compact complex manifold with $\mbox{dim}_\C X=n$. For every $r\in\N^\star$ and every $k=0,\dots , 2n$, the following identity (a kind of numerical Poincar\'e duality extended to all the pages of the spectral sequence) holds:

$$\mbox{dim}_\C E_r^k = \mbox{dim}_\C E_r^{2n-k},$$ 

\noindent where, as usual, $E_r^k=\sum_{p+q=k}E_r^{p,\,q}$ is the direct sum of the spaces of total degree $k$ on the $r^{th}$ page of the Fr\"olicher spectral sequence of $X$.

\end{Prop}

This is an immediate consequence of Theorem \ref{The:main1} and of the following

\begin{Prop}\label{Prop:eigenspace-duality} Let $(X,\,\omega)$ be a compact complex Hermitian manifold with $\mbox{dim}_\C X=n$. Fix an arbitrary constant $h>0$.

\vspace{1ex}

$(i)$\, If $d_h^\star$, resp. $\star$, are the formal adjoint of $d_h$, resp. the Hodge star operator induced by $\omega$, then $$d_h^\star = -\star\bar{d}_h\star.$$

$(ii)$\, If, for every $h>0$, every $k=0,\dots , 2n$ and every $\lambda\geq 0$, $E^k_{\Delta_h}(\lambda)$ stands for the $\lambda$-eigenspace of $\Delta_h:C^\infty_k(X,\,\C)\longrightarrow C^\infty_k(X,\,\C)$, the linear map

$$E^k_{\Delta_h}(\lambda)\longrightarrow E^{2n-k}_{\Delta_h}(\lambda), \hspace{2ex} u\mapsto\star\bar{u},$$

\noindent is well defined and an isomorphism.

 In particular, the operators $\Delta_h:C^\infty_k(X,\,\C)\longrightarrow C^\infty_k(X,\,\C)$ and $\Delta_h:C^\infty_{2n-k}(X,\,\C)\longrightarrow C^\infty_{2n-k}(X,\,\C)$ have the same spectra and their corresponding eigenvalues have the same multiplicities for all $h>0$ and all $k=0,\dots , 2n$.

\end{Prop}

\noindent {\it Proof.} $(i)$\, We have $d_h^\star = h\partial^\star + \bar\partial^\star = -h\star\bar\partial\star - \star\partial\star = -\star(h\bar\partial + \partial)\star = -\star\bar{d}_h\star$ thanks to the standard formulae $\partial^\star = -\star\bar\partial\star$ and $\bar\partial^\star = -\star\partial\star$. 

$(ii)$\, Using the formula under $(i)$ and $\star\star = (-1)^k$ on $k$-forms, we get the following equivalences: \begin{eqnarray}\nonumber u\in E^k_{\Delta_h}(\lambda) & \iff & -d_h\star\bar{d}_h\star u - \star\bar{d}_h\star d_h u = \lambda u \\
\nonumber & \stackrel{(a)}{\iff} & (-\star \bar{d}_h\star)d_h(\star\bar{u}) - (-1)^{\deg u}\,\star\star d_h\star\bar{d}_h\star\star\bar{u} = \lambda\,(\star\bar{u}) \\
\nonumber & \iff & d_h^\star d_h(\star\bar{u}) + d_h d_h^\star(\star\bar{u}) = \lambda\,(\star\bar{u}) \iff \star\bar{u}\in E^{2n-k}_{\Delta_h}(\lambda),\end{eqnarray}

\noindent where (a) was obtained by conjugating and then applying the isomorphism $\star$.

 This shows the well-definedness of the linear map under consideration. Both the conjugation and $\star$ are isomorphisms, hence so is that linear map.  \hfill $\Box$

\vspace{3ex}

\noindent {\it Proof of Proposition \ref{Prop:numerical-F-duality}.} By Theorem \ref{The:main1}, $\mbox{dim}_\C E_r^k$, resp. $\mbox{dim}_\C E_r^{2n-k}$, is the number of eigenvalues $\lambda_i^k(h)\in{\cal O}(h^{2r})$, resp. $\lambda_i^{2n-k}(h)\in{\cal O}(h^{2r})$, counted with multiplicities, of $\Delta_h$ in degree $k$, resp. $2n-k$. Since, by Proposition \ref{Prop:eigenspace-duality}, $\lambda_i^k(h) = \lambda_i^{2n-k}(h)$ for all $i\in\N^\star$ and all $h>0$, the statement follows.  \hfill $\Box$

\vspace{3ex}

The last consequence of Theorem \ref{The:main1} that we notice in this section is the following degeneration criterion for the Fr\"olicher spectral sequence. 

\begin{Prop}\label{Prop:degeneration-criterion} Let $(X,\,\omega)$ be a compact complex Hermitian manifold with $\mbox{dim}_\C X=n$. For every constant $h>0$, let $\delta_h^{(k)}>0$ be the smallest positive eigenvalue of $\Delta_h:C^\infty_k(X,\,\C)\longrightarrow C^\infty_k(X,\,\C)$.

Then, for every $r\in\N^\star$, the Fr\"olicher spectral sequence of $X$ degenerates at $E_r$ if and only if $$\lim\limits_{h\rightarrow 0}\frac{\delta_h^{(k)}}{h^{2r}} = +\infty, \hspace{3ex} \mbox{for all}\hspace{2ex} k\in\{1,\dots , n\}.$$

\end{Prop}

\noindent {\it Proof.} The multiplicity of $0$ as an eigenvalue of $\Delta_h:C^\infty_k(X,\,\C)\longrightarrow C^\infty_k(X,\,\C)$ is the $k^{th}$ Betti number $b_k$ of $X$ (cf. Corollary \ref{Cor:Hodge_d_h}), so the degeneration at $E_r$ of the Fr\"olicher spectral sequence (known to be equivalent to the identities $b_k=\mbox{dim}\,E_r^k$ for all $k=0,1,\dots , 2n$) amounts, thanks to Theorem \ref{The:main1}, to $\delta_h^{(k)}$ converging to zero (if it does converge to zero at all as $h\downarrow 0$) strictly less fast than $Ch^{2r}$ for all $k=0,1,\dots , 2n$. On the other hand, the numerical duality statement of Proposition \ref{Prop:numerical-F-duality} reduces the verification of this property to the cases $k=1,\dots , n$.   \hfill $\Box$

\section{Degeneration at $E_2$ of the Fr\"olicher spectral sequence}\label{section:degeneration_E2}

In this section, we prove Theorem \ref{The:main2}.

We start off by noticing a lower estimate for $\Delta_h-h^2\Delta$ that holds for any Hermitian metric.

\begin{Lem}\label{Lem:Delta_h_Delta_comparison_SKT_3} Let $(X,\,\omega)$ be a compact complex manifold. For every $0<h<1$, the following inequality of operators holds on smooth differential forms of all degrees:

\begin{equation}\label{eqn:Delta_h_Delta_comparison_SKT_3}\Delta_h-h^2\Delta\geq (1-h)h\,\bigg(\Delta''-h[\tau,\,\tau^\star]\bigg).\end{equation}

\end{Lem}

\noindent {\it Proof.} We know from Lemma \ref{Lem:relations_Laplacians} that $\Delta_h = h^2\Delta' + \Delta'' - h[\tau,\,\bar\partial^\star] - h[\bar\partial,\,\tau^\star]$ for any Hermitian metric $\omega$, while $\Delta=[\partial+\bar\partial,\,\partial^\star+\bar\partial^\star] = \Delta' + \Delta'' -[\tau,\,\bar\partial^\star] - [\bar\partial,\,\tau^\star]$. Thus, we get \begin{eqnarray}\label{eqn:Delta_h_Delta_difference}\nonumber\Delta_h - h^2\Delta & = & (1-h^2)\,\Delta'' + h(h-1)\,([\bar\partial,\,\tau^\star] + [\bar\partial^\star,\,\tau]) \\
 & = & (1-h)\,\bigg((1+h)\,\Delta'' -h\,[\bar\partial,\,\tau^\star] -h\,[\bar\partial^\star,\,\tau]\bigg).\end{eqnarray} We shall now estimate the signless terms on the r.h.s. of (\ref{eqn:Delta_h_Delta_difference}). For any form $u$, we have \begin{eqnarray}\nonumber\langle\langle[\bar\partial,\,\tau^\star]\,u,\,u\rangle\rangle + \langle\langle[\bar\partial^\star,\,\tau]\,u,\,u\rangle\rangle & = & \langle\langle\tau^\star u,\,\bar\partial^\star u\rangle\rangle + \langle\langle\bar\partial u,\,\tau u\rangle\rangle + \langle\langle\tau u,\,\bar\partial u\rangle\rangle + \langle\langle\bar\partial^\star u,\,\tau^\star u\rangle\rangle\\
\nonumber & = & 2\mbox{Re}\,\langle\langle\bar\partial^\star u,\,\tau^\star u\rangle\rangle + 2\mbox{Re}\,\langle\langle\bar\partial u,\,\tau u\rangle\rangle.\end{eqnarray}

\noindent Thus, for any Hermitian metric $\omega$, we have \begin{eqnarray}\nonumber h\,|\langle\langle([\bar\partial,\,\tau^\star] + [\bar\partial^\star,\,\tau])\,u,\,u\rangle\rangle| & \leq & 2h\,|\langle\langle\bar\partial u,\,\tau u\rangle\rangle| + 2h\,|\langle\langle\bar\partial^\star u,\,\tau^\star u\rangle\rangle|\\
\nonumber & \leq & (||\bar\partial u||^2 + ||\bar\partial^\star u||^2) + h^2\,(||\tau u||^2 + ||\tau^\star u||^2 ) \\
\nonumber & = & \langle\langle\Delta''u,\,u\rangle\rangle + h^2\,\langle\langle[\tau,\,\tau^\star]\,u,\,u\rangle\rangle.\end{eqnarray}

 Using this last estimate in (\ref{eqn:Delta_h_Delta_difference}), we get $\Delta_h - h^2\Delta \geq (1-h)\,(h\Delta'' - h^2[\tau,\,\tau^\star])$ in the sense of operators. This is precisely (\ref{eqn:Delta_h_Delta_comparison_SKT_3}).   

 Note that we can also write $|\langle\langle([\bar\partial,\,\tau^\star] + [\bar\partial^\star,\,\tau])\,u,\,u\rangle\rangle|\leq \langle\langle\Delta''u,\,u\rangle\rangle + \langle\langle[\tau,\,\tau^\star]\,u,\,u\rangle\rangle$ for every form $u$, which, alongside (\ref{eqn:Delta_h_Delta_difference}), yields $\Delta_h - h^2\Delta \geq (1-h)\,(\Delta'' - h[\tau,\,\tau^\star])$. This is slightly better than (\ref{eqn:Delta_h_Delta_comparison_SKT_3}) if the r.h.s. is non-negative, but worse otherwise. \hfill $\Box$

\vspace{3ex}

We shall now give a sufficient condition for the r.h.s. of (\ref{eqn:Delta_h_Delta_comparison_SKT_3}) to be non-negative.

\begin{Lem}\label{Lem:ker-inclusion_l-bound} Let $(X,\,\omega)$ be a compact Hermitian manifold with $\mbox{dim}_\C X=n$ such that the inclusion of kernels $$\ker\Delta''\subset\ker\,[\tau,\,\tau^\star]$$

\noindent holds for the operators $\Delta'', [\tau,\,\tau^\star]:C^\infty_k(X,\,\C)\longrightarrow C^\infty_k(X,\,\C)$ in a fixed degree $k\in\{1,\dots , n\}$. 

Then, there exists a constant $h_0(k)\in(0,\,1]$ such that the following inequality of operators holds in degree $k$: $$\Delta''\geq h\,[\tau,\,\tau^\star]  \hspace{3ex} \mbox{for all}\hspace{1ex} 0<h<h_0(k).$$

\end{Lem}

\noindent {\it Proof.} Let $\delta''_k>0$ be the smallest positive eigenvalue of the elliptic, self-adjoint and non-negative differential operator $\Delta'':C^\infty_k(X,\,\C)\longrightarrow C^\infty_k(X,\,\C)$. 

 On the other hand, the operator $[\tau,\,\tau^\star]:C^\infty_k(X,\,\C)\longrightarrow C^\infty_k(X,\,\C)$ is of order zero, hence bounded, so the constant $C_k:=\sup_{||u||\leq 1}\langle\langle[\tau,\,\tau^\star]\,u,\,u\rangle\rangle$ is finite. 

We put $h_0(k):=\min\{\delta''_k/C_k,\,1\}$ and will prove that $\langle\langle\Delta''u,\,u\rangle\rangle\geq h\,\langle\langle[\tau,\,\tau^\star]\,u,\,u\rangle\rangle$ for all $u\in C^\infty_k(X,\,\C)$ and all $h\in(0,\,h_0(k))$. Let us fix a form $u\in C^\infty_k(X,\,\C)$. 

 Since $\Delta''$ is elliptic and preserves bidegrees, the following orthogonal splitting $$C^\infty_k(X,\,\C) = \ker\Delta''\oplus\mbox{Im}\,\Delta''$$ 

\noindent holds and induces a unique splitting $u=u_h+u_{h^\perp}$ with $u_h\in\ker\Delta''$ and $u_{h^\perp}\in\mbox{Im}\,\Delta''$. In particular, $u_h\in\ker\,[\tau,\,\tau^\star]$ thanks to our assumption.  

 We get \begin{equation}\label{eqn:l-bound_Delta''}\langle\langle\Delta''u,\,u\rangle\rangle = \langle\langle\Delta''u_{h^\perp},\,u_h+u_{h^\perp}\rangle\rangle = \langle\langle\Delta''u_{h^\perp},\,u_{h^\perp}\rangle\rangle\geq\delta''_k\,||u_{h^\perp}||^2\end{equation}

\noindent since $u_{h^\perp}\perp\ker\Delta''$, so $u_{h^\perp}$ lies in the orthogonal direct sum of the eigenspaces of $\Delta''$ corresponding to positive eigenvalues ($=$ eigenvalues $\geq\delta''_k$). 

 On the other hand, \begin{eqnarray}\label{eqn:u-bound_torsion}\nonumber\langle\langle[\tau,\,\tau^\star]\,u,\,u\rangle\rangle & \stackrel{(a)}{=} & \langle\langle[\tau,\,\tau^\star]\,u_{h^\perp},\,u_h + u_{h^\perp}\rangle\rangle \stackrel{(b)}{=} \langle\langle u_{h^\perp},\,[\tau,\,\tau^\star]\,u_h\rangle\rangle + \langle\langle[\tau,\,\tau^\star]\,u_{h^\perp},\,u_{h^\perp}\rangle\rangle \\
 & \stackrel{(c)}{=} & \langle\langle[\tau,\,\tau^\star]\,u_{h^\perp},\,u_{h^\perp}\rangle\rangle \stackrel{(d)}{\leq} C_k\, ||u_{h^\perp}||^2,\end{eqnarray}

\noindent where for (a) we used the fact that $u_h\in\ker\,[\tau,\,\tau^\star]$, for (b) we used the self-adjointness of $[\tau,\,\tau^\star]$, (c) follows from $u_h\in\ker\,[\tau,\,\tau^\star]$, while (d) follows from the definition of $C_k$. 

 Since $h_0(k)=\min\{\frac{\delta''_k}{C_k},\,1\}$, inequalities (\ref{eqn:l-bound_Delta''}) and (\ref {eqn:u-bound_torsion}) imply that

$$h\,\langle\langle [\tau,\,\tau^\star]\,u,\,u\rangle\rangle\leq C_kh\, ||u_{h^\perp}||^2 \leq\frac{C_kh}{\delta''_k}\,\langle\langle\Delta''u,\,u\rangle\rangle\leq\langle\langle\Delta''u,\,u\rangle\rangle$$ 

\noindent for all $h\in(0,\,h_0(k))$.   \hfill $\Box$

\begin{Cor}\label{Cor:Delta_h_h2_Delta} Let $(X,\,\omega)$ be a compact Hermitian manifold such that $\ker\Delta''\subset\ker\,[\tau,\,\tau^\star]$ in a fixed degree $k$. Then, there exists a constant $h_0(k)\in(0,\,1]$ such that the following inequality of operators holds in degree $k$: $$\Delta_h\geq h^2\,\Delta  \hspace{3ex} \mbox{for all}\hspace{1ex} 0<h<h_0(k).$$

\end{Cor}

\noindent {\it Proof.} This is an immediate consequence of Lemmas \ref{Lem:Delta_h_Delta_comparison_SKT_3} and \ref{Lem:ker-inclusion_l-bound}.  \hfill $\Box$

\vspace{3ex}

We can now prove the spectral sequence degeneration statement of this paper.

\vspace{3ex}

\noindent {\bf Proof of Theorem \ref{The:main2}.} Let us fix an arbitrary $k\in\{1,\dots , n\}$. Hypothesis (\ref{eqn:torsion_hyp}) and Corollary \ref{Cor:Delta_h_h2_Delta} imply that $\ker\Delta_h\subset\ker\Delta$ for all $0<h<h_0(k)$ since $\langle\langle\Delta u,\,u\rangle\rangle\geq 0$ for every $u$ and $u\in\ker\Delta$ if and only if $\langle\langle\Delta u,\,u\rangle\rangle= 0$. Meanwhile, we know from Corollary \ref{Cor:Hodge_d_h} that $\ker\Delta_h$ and $\ker\Delta$ are finite-dimensional vector spaces of {\bf equal dimensions}, so for all $0<h<h_0(k)$ we get

\begin{equation}\label{eqn:kernels-equality}\ker\Delta_h = \ker\Delta.\end{equation}

 For every $h>0$, let $\delta_h^{(k)}>0$ be the smallest positive eigenvalue of the elliptic operator $\Delta_h:C^\infty_k(X,\,\C)\longrightarrow C^\infty_k(X,\,\C)$ and let $u_h\in C^\infty_k(X,\,\C)$ be a corresponding unitary eigenvector, i.e. $$||u_h||=1 \hspace{3ex} \mbox{and} \hspace{3ex} \Delta_hu_h=\delta_h^{(k)}u_h.$$

\noindent Now, $u_h$ is orthogonal to $\ker\Delta_h$, hence, thanks to (\ref{eqn:kernels-equality}), $u_h$ is also orthogonal to $\ker\Delta$ for every $0<h<h_0(k)$. Consequently, $\langle\langle\Delta u_h,\,u_h\rangle\rangle\geq\delta_k\,||u_h||^2 = \delta_k$, where $\delta_k>0$ is the smallest positive eigenvalue of $\Delta:C^\infty_k(X,\,\C)\longrightarrow C^\infty_k(X,\,\C)$.

Using this and Corollary \ref{Cor:Delta_h_h2_Delta}, we get $$\delta_h^{(k)}=\langle\langle\Delta_hu_h,\,u_h\rangle\rangle \geq h^2\,\langle\langle\Delta u_h,\,u_h\rangle\rangle\geq\delta_k h^2 \hspace{3ex} \mbox{for all}\hspace{1ex} 0<h<h_0(k).$$

\noindent In particular, $\lim_{h\rightarrow 0}(\delta_h^{(k)}/h^4) = +\infty$. 

As in the proof of Proposition \ref{Prop:degeneration-criterion}, this and Theorem \ref{The:main1} imply that $\mbox{dim}\,E_2^k = b_k$ for the degree $k\in\{1,\dots , n\}$ that was arbitrarily fixed in the beginning. By the duality statement of Proposition \ref{Prop:numerical-F-duality}, this also yields $\mbox{dim}\,E_2^{2n-k} = b_k = b_{2n-k}$. Since this holds for all $k\in\{1,\dots , n\}$, the Fr\"olicher spectral sequence of $X$ degenerates at $E_2$.  \hfill $\Box$

\section{Appendix: Comparison of Laplacians when the metric is SKT}\label{section:appendix}

In this section, we come within an $\varepsilon$ ($=Ch^2$) of solving Conjecture \ref{Conj:SKT_E2} as an application of Theorem \ref{The:main1} and of a comparison of the Laplacians $\Delta'$ and $\Delta''$ defined by an arbitrary SKT metric $\omega$ supposed to exist on a given compact complex manifold $X$. Recall that an SKT metric $\omega$ is a $C^\infty$ positive definite $(1,\,1)$-form $\omega$ such that $\partial\bar\partial\omega=0$ on $X$.

\begin{Lem}\label{Lem:Delta'-''_comparison_SKT} Let $X$ be a compact complex manifold on which an {\bf SKT metric $\omega$} exists.

\vspace{1ex}

 $(i)$\, The usual $\partial$- and $\bar\partial$-Laplacians $\Delta'=[\partial,\,\partial^\star]$ and $\Delta''=[\bar\partial,\,\bar\partial^\star]$ induced by $\omega$ satisfy the following inequalities on differential forms of all bidegrees: \begin{eqnarray}\label{eqn:Delta'-''_comparison_SKT_simplified}(1+\delta)\,\Delta'' + \bigg(1+\frac{1}{\delta}\bigg)\,[\bar\tau,\,\bar\tau^\star]\geq\Delta'\geq\frac{1}{1+\delta}\,\Delta'' - \frac{1}{\delta}\,[\tau,\,\tau^\star], \hspace{3ex} \mbox{for all} \hspace{1ex} \delta>0,\end{eqnarray} \noindent where $\tau=\tau_\omega:=[\Lambda_\omega,\,\partial\omega\wedge\cdot]$ is the torsion operator of type $(1,\,0)$ and $\bar\tau^\star$ is the formal adjoint w.r.t. the $L^2_\omega$-inner product of its complex conjugate.

\vspace{1ex}

$(ii)$\, The following inequality also holds: \begin{eqnarray}\label{eqn:Delta'-''_comparison_SKT_h}\Delta''\geq h\,\Delta' + \bigg(h\overline{X}_\omega - \frac{h}{1-h}\,[\bar\tau,\,\bar\tau^\star]\bigg), \hspace{3ex} \mbox{for all} \hspace{1ex} 0<h<1,\end{eqnarray} \noindent where $X_\omega:=[\partial\omega\wedge\cdot,\,(\partial\omega\wedge\cdot)^\star]$. Implicitly, we have \begin{eqnarray}\label{eqn:Delta_h_Delta_comparison_SKT}\Delta_h-h\Delta\geq h\,\bigg((1-h)\,\overline{X}_\omega - [\bar\tau,\,\bar\tau^\star]\bigg), \hspace{3ex} \mbox{for all} \hspace{1ex} 0<h<1.\end{eqnarray}

\noindent Since $\overline{X}_\omega$ and $[\bar\tau,\,\bar\tau^\star]$ are zero-th order operators, they are bounded, so (\ref{eqn:Delta_h_Delta_comparison_SKT}) implies the existence of a constant $C>0$ independent of $h$ such that \begin{eqnarray}\label{eqn:Delta_h_Delta_comparison_SKT_C}\Delta_h-h\Delta\geq -Ch, \hspace{3ex} \mbox{for all} \hspace{1ex} 0<h<1.\end{eqnarray}

\end{Lem}

\noindent {\it Proof.} $(i)$\, Demailly's formula (cf. [Dem84] or [Dem97, VII, $\S.1$]) of the Bochner-Kodaira-Nakano type for arbitrary Hermitian metrics $\omega$ reads

$$\Delta' = \Delta''_{\bar\tau} - \overline{X}_\omega + [\Lambda_\omega,\,[\Lambda_\omega,\,\frac{i}{2}\,\partial\bar\partial\omega]],$$

\noindent where $\Delta''_{\bar\tau}:=[\bar\partial+\bar\tau,\,(\bar\partial+\bar\tau)^\star]$ and $\overline{X}_\omega:=[\bar\partial\omega\wedge\cdot,\,(\bar\partial\omega\wedge\cdot)^\star]$. The last term on the r.h.s. above vanishes if $\omega$ is SKT, so we get

\begin{equation}\label{eqn:BKN_1}\Delta'' + ([\bar\partial,\,\bar\tau^\star] + [\bar\tau,\,\bar\partial^\star]) + [\bar\tau,\,\bar\tau^\star] = \Delta' + \overline{X}_\omega  \hspace{5ex} \mbox{if}\hspace{1ex} \partial\bar\partial\omega=0.\end{equation}

 Now, the signless terms can be easily estimated using the elementary inequality $2|ab|\leq \delta a^2 + (1/\delta)\,b^2$ for arbitrary $a,b\in\C$ and $\delta>0$. For every differential form $u$ of any degree, we get: 

\begin{eqnarray}\nonumber & & |\langle\langle[\bar\partial,\,\bar\tau^\star]\,u,\,u\rangle\rangle + \langle\langle[\bar\tau,\,\bar\partial^\star]\,u,\,u\rangle\rangle| = |2\mbox{Re}\,\langle\langle\bar\partial u,\,\bar\tau u\rangle\rangle + 2\mbox{Re}\,\langle\langle\bar\partial^\star u,\,\bar\tau^\star u\rangle\rangle| \\
\nonumber & \leq & 2\,|\langle\langle\bar\partial u,\,\bar\tau u\rangle\rangle| + 2\,|\langle\langle\bar\partial^\star u,\,\bar\tau^\star u\rangle\rangle| \leq \delta\,||\bar\partial u||^2 + \frac{1}{\delta}\,||\bar\tau u||^2 + \delta\,||\bar\partial^\star u||^2 + \frac{1}{\delta}\,||\bar\tau^\star u||^2  \\
\nonumber & = & \delta\,\langle\langle\Delta'' u,\,u\rangle\rangle + \frac{1}{\delta}\,\langle\langle[\bar\tau,\,\bar\tau^\star]\, u,\,u\rangle\rangle.\end{eqnarray}

\noindent Together with (\ref{eqn:BKN_1}), this implies that $(1+\delta)\,\Delta'' + (1+ 1/\delta)\,[\bar\tau,\,\bar\tau^\star]\geq\Delta' + \overline{X}_\omega$ if $\omega$ is SKT. This is essentially an upper estimate for $\Delta'$ whose conjugate yields a lower estimate for $\Delta'=\overline{\Delta''}$. Putting these upper and lower estimates together, we get

\begin{equation}\label{eqn:Delta'-''_comparison_SKT}(1+\delta)\,\Delta'' + \bigg(1+\frac{1}{\delta}\bigg)\,[\bar\tau,\,\bar\tau^\star] - \overline{X}_\omega \geq\Delta'\geq\frac{1}{1+\delta}\,\Delta'' + \frac{1}{1+\delta}\, X_\omega - \frac{1}{\delta}\,[\tau,\,\tau^\star],\end{equation}

\noindent for all $\delta>0$. Since $X_\omega$ and $\overline{X}_\omega$ are non-negative operators, ignoring them weakens these inequalities to (\ref{eqn:Delta'-''_comparison_SKT_simplified}).

\vspace{2ex}

$(ii)$\, After dividing by $1+\delta$, the l.h.s. inequality in (\ref{eqn:Delta'-''_comparison_SKT}) translates to

$$\Delta'' \geq \frac{1}{1+\delta}\,\Delta' + \frac{1}{1+\delta}\,\overline{X}_\omega - \frac{1}{\delta}\,[\bar\tau,\,\bar\tau^\star].$$

\noindent This is precisely (\ref{eqn:Delta'-''_comparison_SKT_h}) if we put $h:=\frac{1}{1+\delta}\in(0,\,1)$ since in this case $\delta=\frac{1-h}{h}$.

 To get (\ref{eqn:Delta_h_Delta_comparison_SKT}) from (\ref{eqn:Delta'-''_comparison_SKT_h}), it suffices to notice that $\Delta_h-h\Delta = h(h-1)\,\Delta' + (1-h)\,\Delta'' = (1-h)\,(\Delta'' - h\Delta')$.  \hfill $\Box$

\vspace{3ex}

We now observe an analogue of inequality (\ref{eqn:Delta_h_Delta_comparison_SKT}) for $\Delta_h-h^2\Delta$. 

\begin{Lem}\label{Lem:Delta'-''_comparison_SKT_2} Let $X$ be a compact complex manifold on which an {\bf SKT metric $\omega$} exists. The following inequalities of operators hold: \begin{eqnarray}\label{eqn:Delta_h_Delta_comparison_SKT_2}\Delta_h-h^2\Delta\geq h^2\,\bigg((1-h)\,\overline{X}_\omega - [\bar\tau,\,\bar\tau^\star]\bigg) \geq -Ch^2, \hspace{3ex} \mbox{for all} \hspace{1ex} 0<h<1,\end{eqnarray} \noindent where $\overline{X}_\omega:=[\bar\partial\omega\wedge\cdot,\,(\bar\partial\omega\wedge\cdot)^\star]$ and $C\geq 0$ is a constant independent of $h$.

\end{Lem}

\noindent {\it Proof.} Since $\Delta_h = h^2\Delta' + \Delta'' + hA$ and $\Delta = \Delta' + \Delta'' + A$, where $A:=[\partial,\,\bar\partial^\star] + [\bar\partial,\,\partial^\star]$, we get

$$\Delta_h - h^2\Delta = (1-h)\,((1+h)\,\Delta'' + hA).$$ 

 On the other hand, the signless operator $A$ can be estimated in the same way as a similar operator was estimated in the proof of Lemma \ref{Lem:Delta'-''_comparison_SKT}. We get $\langle\langle Au,\,u\rangle\rangle = 2\mbox{Re}\,\langle\langle\partial u,\,\bar\partial u\rangle\rangle + 2\mbox{Re}\,\langle\langle\partial^\star u,\,\bar\partial^\star u\rangle\rangle$, hence $$h\,|\langle\langle Au,\,u\rangle\rangle| \leq h^2\,||\partial u||^2 + ||\bar\partial u||^2 + h^2\,||\partial^\star u||^2 + ||\bar\partial^\star u||^2 = h^2\,\langle\langle\Delta'u,\,u\rangle\rangle + \langle\langle\Delta''u,\,u\rangle\rangle$$ \noindent for any form $u$. Consequently, $(1+h)\,\Delta'' + hA\geq h\Delta'' - h^2\Delta'$ as operators, so we get

$$\Delta_h - h^2\Delta\geq h(1-h)\,(\Delta'' - h\Delta').$$

\noindent (Note that we can also write $|\langle\langle Au,\,u\rangle\rangle| \leq \langle\langle\Delta'u,\,u\rangle\rangle + \langle\langle\Delta''u,\,u\rangle\rangle$ and we get $\Delta_h - h^2\Delta = (1-h)\,((1+h)\,\Delta'' + hA) \geq (1-h)\,(\Delta'' - h\Delta')$ for every form $u$.)

 Meanwhile, from (\ref{eqn:Delta'-''_comparison_SKT_h}) we know that $(1-h)\,(\Delta''- h\,\Delta') \geq h\,\bigg((1-h)\,\overline{X}_\omega - [\bar\tau,\,\bar\tau^\star]\bigg)$ for all $0<h<1$. Together with the last inequality, this proves the first inequality in (\ref{eqn:Delta_h_Delta_comparison_SKT_2}). 

 The second inequality in (\ref{eqn:Delta_h_Delta_comparison_SKT_2}) follows at once from the first since $\overline{X}_\omega\geq 0$ and the non-negative operator $[\bar\tau,\,\bar\tau^\star]$ is of order zero, hence bounded, so we can choose $C:=\sup_{||u||=1}\langle\langle[\bar\tau,\,\bar\tau^\star]u,\,u\rangle\rangle<+\infty$. 

(Using the alternative lower estimate $\Delta_h - h^2\Delta \geq (1-h)\,(\Delta'' - h\Delta')$ noticed above, the inequalities in (\ref{eqn:Delta_h_Delta_comparison_SKT_2}) get replaced by $\Delta_h-h^2\Delta\geq h\,((1-h)\,\overline{X}_\omega - [\bar\tau,\,\bar\tau^\star])\geq -Ch$.)

\hfill $\Box$

\vspace{3ex}

If the lower bound $-Ch^2$ in (\ref{eqn:Delta_h_Delta_comparison_SKT_2}) could be improved to $0$, then we would have $\Delta_h\geq h^2\Delta$ for all $0<h\ll 1$ (as in Corollary \ref{Cor:Delta_h_h2_Delta}) and Conjecture \ref{Conj:SKT_E2} would follow by the argument spelt out at the end of section $\S.$\ref{section:degeneration_E2}.

\vspace{3ex}

\noindent {\bf References.} \\

\noindent [ALK00]\, J.A. \'Alvarez L\'opez, Y.A. Kordyukov --- {\it Adiabatic Limits and Spectral Sequences for Riemannian Foliations} ---  Geom. Funct. Anal. 10 (2000), no. 5, 977–1027. 

\vspace{1ex}

\noindent [CFGU97]\, L.A. Cordero, M. Fernandez, A. Gray, L. Ugarte --- {\it A General Description of the Terms in the Fr\"olicher Spectral Sequence} --- Diff. Geom. and its Applic. {\bf 7} (1997), 75–84.

\vspace{1ex}

\noindent [Dem 84]\, J.-P. Demailly --- {\it Sur l'identit\'e de Bochner-Kodaira-Nakano en g\'eom\'etrie hermitienne} --- S\'eminaire d'analyse P. Lelong, P. Dolbeault, H. Skoda (editors) 1983/1984, Lecture Notes in Math., no. {\bf 1198}, Springer Verlag (1986), 88-97.

\vspace{1ex}

\noindent [Dem 97]\, J.-P. Demailly --- {\it Complex Analytic and Algebraic Geometry}---http://www-fourier.ujf-grenoble.fr/~demailly/books.html

\vspace{1ex}

\vspace{1ex}

\noindent [ES89]\, D.V. Efremov, M.A. Shubin --- {\it Spectrum Distribution Function and Variational Principle for Automorphic Operators on Hyperbolic Space} --- 

\vspace{1ex}

\noindent [For95]\, R. Forman --- {\it Spectral Sequences and Adiabatic Limits} --- Commun. Math. Phys. {\bf 168}, 57-116 (1995).

\vspace{1ex}

\noindent [Fro55]\, A. Fr\"olicher --- {\it Relations between the Cohomology Groups of Dolbeault and Topological Invariants} --- Proc. Nat. Acad. Sci. U.S.A. 41 (1955), 641–644. 

\vspace{1ex}

\noindent [GS91]\, M. Gromov, M.A. Shubin --- {\it Von Neumann Spectra near Zero} ---  Geom. Funct. Anal. 1 (1991), no. 4, 375–404.

\vspace{1ex}

\noindent [MM90]\, R. R. Mazzeo, R. B. Melrose --- {\it The Adiabatic Limit, Hodge Cohomology and Leray's Spectral Sequence} --- J. Diff. Geom. {\bf 31} (1990) 185-213.

\vspace{1ex}

\noindent [Pop16]\, D. Popovici --- {\it Degeneration at $E_2$ of Certain Spectral Sequences} ---  International Journal of Mathematics {\bf 27}, no. 14 (2016), DOI: 10.1142/S0129167X16501111.

\vspace{1ex}

\noindent [Wi85]\, E. Witten --- {\it Global Gravitational Anomalies} --- Commun. Math. Phys, {\bf 100}, 197-229 (1985).

\vspace{6ex}

\noindent Universit\'e Paul Sabatier, Institut de Math\'ematiques de Toulouse,

\noindent 118 route de Narbonne, 31062 Toulouse, France

\noindent Email: popovici@math.univ-toulouse.fr

\end{document}